\documentclass[12pt]{article}
\usepackage{amsthm, amsmath,amssymb,amsfonts}
\usepackage{graphicx}
\usepackage{color}
\usepackage{psfig}
\usepackage{graphics}
\usepackage{amsmath}
\usepackage{amssymb}
\usepackage{epsfig}
\usepackage{epstopdf}
\usepackage{color}

\newtheorem{theorem}{Theorem}[section]
\newtheorem{proposition}[theorem]{Proposition}
\newtheorem{corollary}[theorem]{Corollary}
\newtheorem{lemma}[theorem]{Lemma}
\newtheorem{remark}[theorem]{Remark}
\newtheorem{definition}[theorem]{Definition}
\newtheorem{example}[theorem]{Example}

\def\para{\vspace{3mm}}

\def\gcd{{\rm gcd}}

\def\Resultant{{\rm Res}}

\def\deg{{\rm deg}}
\def\degree{{\rm deg}}

\def\resultant{{\rm resultant}}

\begin{document}
\title{Asymptotes and Perfect Curves}
\author{Angel Blasco and Sonia P\'erez-D\'{\i}az\\
Departamento de F\'{\i}sica y Matem\'aticas \\
        Universidad de Alcal\'a \\
      E-28871 Madrid, Spain  \\
angel.blasco@uah.es, sonia.perez@uah.es
}
\date{}          % Enter your date or \today between curly braces
\maketitle

\begin{abstract}
We develop a method for computing all the {\it generalized
asymptotes} of a real plane algebraic curve $\cal C$ over $\Bbb C$
implicitly  defined by an irreducible polynomial $f(x,y)\in {\Bbb R}[x,y]$. The approach is based on the notion of perfect curve introduced from
the concepts and results presented in \cite{paper1}.
\end{abstract}


 \section{Introduction}

Let $\cal C$ be a real plane algebraic curve over $\Bbb C$ implicitly defined
 by an irreducible polynomial $f(x,y)\in {\Bbb R}[x,y]$. In this paper, we deal with the problem of computing  the asymptotes of the infinity branches of $\cal C$. This question is  very important in the study of real plane algebraic
curves because  asymptotes contain much of the information about the behavior of the curves in the large. For instance, determining the asymptotes of a curve is an important step in sketching its graph.

\para

Intuitively speaking, the asymptotes of some branch of a real plane
algebraic curve reflect the status of this branch at the points with
sufficiently large coordinates. In analytic geometry, an asymptote
of a curve is a line such that the distance between the curve and
the line approaches zero as they tend to infinity.  In some
contexts, such as algebraic geometry, an asymptote is defined as a
line which is tangent to a curve at infinity.

\para

More precisely, let $\cal C$  be a real plane algebraic curve, and $B$ an infinity branch of $\cal C$. A line $\ell$ is called the asymptote of $\cal C$ at $B$, if for every $\epsilon\in {\Bbb R}^+$, there exists $M\in {\Bbb R}^+$ such that $d(P,\ell)<\epsilon$, for every $P\in B$ with $\|P\|>M.$

\para

If $B$ can be defined by some explicit equation of the form
$y = f(x)$ (or $x = g(y)$), where $f$ (or $g$) is a continuous function on an infinite interval, it is easy
to decide whether $\cal C$  has an asymptote at $B$ by analyzing the existence  of the limits
of certain functions when $x\rightarrow\infty$ (or $y\rightarrow\infty$). Moreover, if these limits  can be  computed, we may obtain the equation of the asymptote
of $\cal C$ at $B$. However, if this branch $B$ is implicitly defined and its equation cannot be converted into an
explicit form, both the decision and the computation of the asymptote of $\cal C$ at $B$ require some others tools.

%For curves given by the graph of a function $y = f(x)$, there are potentially three kinds of asymptotes: horizontal, vertical and oblique asymptotes. Horizontal asymptotes are horizontal lines that the graph of the function approaches as $x$ tends to $+\infty$ or $-\infty$. Vertical asymptotes are vertical lines near which the function grows without bound. When a linear asymptote is not parallel to the $x$ or $y$-axis, it is called an oblique asymptote.  Asymptotes of these elementary functions can be easily computed.

%\para

%Unlike asymptotes for curves that are graphs of functions, a general curve may have more than two non-vertical asymptotes, and may cross its vertical asymptotes more than once. More generally, one curve is a curvilinear asymptote of another (as opposed to a linear asymptote) if the distance between the two curves tends to zero as they tend to infinity, although usually the term asymptote by itself is reserved for linear asymptotes. The asymptotes of an algebraic curve in the affine plane are the lines that are tangent to the projective curve through a point at infinity. Asymptotes are often considered only for real curves  although they also make sense when defined in this way for curves over an arbitrary field.

\para

Determining the asymptotes of an implicit algebraic plane curve  is a topic considered in many text-books on
analysis (see for instance \cite{maxwell}). In   \cite{jovan}, is presented a fast and a simple method for obtaining the asymptotes of a curve defined by an irreducible polynomial, with emphasis on
second order polynomials.  In \cite{Zeng},  an algorithm for computing all the linear asymptotes of a real plane
algebraic curve  $\cal C$   implicitly defined, is obtained. More precisely,   one may decide whether  a branch of $\cal C$ has
an asymptote, compute all the asymptotes of $\cal C$, and determine those branches whose asymptotes
are the same. By this algorithm, all the asymptotes of  $\cal C$    may be represented
via polynomial real root isolation.

\para

An algebraic plane curve may have more general curves than lines describing the status of a branch at the points
with sufficiently large coordinates. This motivates that in this paper, we are interested in analyzing and computing these {\it  generalized asymptotes}. Intuitively speaking, we say that a curve $\widetilde{{\cal C}}$ is a {\it  generalized asymptote} (or \textit{g-asymptote}) of another curve $\cal C$ if the distance between $\widetilde{{\cal C}}$ and $\cal C$  tends to zero as they tend to infinity, and $\cal C$  can not be approached by a new curve of lower degree.

\para

In this paper,   we present an algorithm for computing all the
g-asymptotes of a real  algebraic plane curve $\cal C$ defined
 by an irreducible polynomial $f(x,y)\in {\Bbb R}[x,y]$
(the assumption of reality is included because of the nature of the
problem, but the theory can be similarly developed for the case of
complex non-real curves).  For this purpose, we use the results in
\cite{paper1}, where  the notions of convergent branches  (that is,
branches that get closer  as they tend to infinity) and approaching
curves are introduced.
%More precisely, in \cite{paper1}, we say that  a  curve $\overline{\cal C}$  {\em
%approaches} ${\cal C}$ at its infinity branch $B$ if the distance between
%$\overline{\cal C}$ and $B$ approaches zero as they tend to infinity.
In addition, in \cite{paper1},  we also provide some results that
characterize   whether two  implicit  algebraic plane
curves approach each other at the infinity, and we present a method
to compare the asymptotic behavior  of two curves (i.e., the behavior
  at the infinity). In particular, we prove that if two
plane curves have the same {\it asymptotic behavior},
the Hausdorff distance between them is finite. Some of these results are summarized in this paper (see Section 2).

\para

The study of approaching curves and convergent branches leads to the
notions of \textit{perfect curve} (a curve of degree $d$ that cannot
be approached by any curve of degree less than $d$) and
\textit{g-asymptote} (a perfect curve that approaches another curve
at an infinity branch). These concepts are introduced in Section 3.
In this section, we also develop an algorithm that computes a
g-asymptote for each infinity branch of a given curve. In Section 4,
we provide some necessary and sufficient conditions for a curve to
be perfect. In particular, we show that a perfect curve admits a
polynomial parametrization. In Section 5, we observe that
``proximity'' is an equivalence relation for the set of perfect
curves. Hence, an infinity branch of a given curve does not have, in
general, a unique g-asymptote, but a whole ``equivalence class''
defined by infinitely many curves. We show that all the
curves in a same class have the same degree $d$, and that any class is
isomorphic to $\mathbb{R}^{d(d-1)/2}$. Finally, we present some results that
allows us to obtain, under certain assumptions, all the curves
within a class.

\section{Notation and Previous Results}\label{S-notation}

In this section, we introduce the notion of {\it infinity branch},  {\it convergent branches}
and  {\it approaching curves},
and we obtain some  properties which allow us to compare  the
behavior of two implicit algebraic plane curves at  the infinity.
For more details on these concepts and results, we refer to
\cite{paper1}.

\para

Throughout the paper,  we consider an algebraic affine plane curve $\cal C$  over
$\mathbb{C}$  defined by the irreducible polynomial $f(x,y) \in {\Bbb R}[x,y]$.
Let ${\cal C}^*$  be its corresponding projective curve  defined by
the homogeneous polynomial $$F(x,y,z)=f_d(x,y)+zf_{d-1}(x,y)+z^2f_{d-2}(x,y)+\cdots+z^d f_0 \in {\Bbb R}[x,y,z],\quad $$
where $d:=\deg({\cal C})$. We assume that $(0:1:0)$ is not an infinity point of  ${\cal C}^*$; otherwise, we may consider a linear change of coordinates.

\para

Let $P=(1:m:0),\,m\in {\Bbb C}$ be an infinity point of  ${\cal C}^*$, and we consider the  curve  defined by the
polynomial $g(y,z)=F(1:y:z)$. We compute the
series expansion for the solutions of $g(y,z)=0$. There exist
exactly $\deg_{Y}(g)$ solutions given by different Puiseux series
that can be grouped into  conjugacy classes.  More precisely, if $$\varphi(z)=m+a_1z^{N_1/N}+a_2z^{N_2/N}+a_3z^{N_3/N} +
\cdots\in{\Bbb C}\ll z\gg,\quad a_i\not=0,\, \forall i\in {\Bbb
N},$$
where $N\in {\Bbb
N}$ , $N_i\in {\Bbb
N},\,\,i=1,\ldots$,  and $0<N_1<N_2<\cdots$, is a Puiseux
series such that $g(\varphi(z), z)=0$, and   $\nu(\varphi)=N$  (i.e., $N$ is the  {\em  ramification index} of $\varphi$), the
 series
 \[\varphi_j(z)=m+a_1c_j^{N_1}z^{N_1/N}+a_2c_j^{N_2}z^{N_2/N}+a_3c_j^{N_3}z^{N_3/N}
+ \cdots\] where $c_j^N=1,\,\,j=1,\ldots,N$,  are called
the {\em conjugates} of $\varphi$.
The set of all (distinct) conjugates of $\varphi$ is called the  {\em  conjugacy
class} of $\varphi$, and the number of different conjugates of $\varphi$ is $\nu(\varphi)$ (see \cite{Duval89}).

\para

Since
$g(\varphi(z), z)=0$ in some neighborhood of $z=0$ where
$\varphi(z)$ converges, there exists $M \in {\Bbb R}^+$
such that
$$F(1:\varphi(t):t)=g(\varphi(t), t)=0,\quad \mbox{for\, $t\in {\Bbb C}$\, and
$|t|<M$},$$
which implies that
$F(t^{-1}:t^{-1}\varphi(t):1)=f(t^{-1},t^{-1}\varphi(t))=0$, for
$t\in {\Bbb C}$  and  $0<|t|<M$. We set $t^{-1}=z$, and   we obtain
that
$$f(z,r(z))=0,\quad \mbox{$z\in {\Bbb C}$\, and
$|z|>M^{-1}$,\qquad where}$$
$$r(z)=z\varphi(z^{-1})=mz+a_1z^{1-N_1/N}+a_2z^{1-N_2/N}+a_3z^{1-N_3/N} + \cdots,\quad a_i\not=0,\, \forall i\in {\Bbb N} $$
$N,N_i\in {\Bbb N},\,\,i=1,\ldots$, and $0<N_1<N_2<\cdots$.

\para

Reasoning similarly with the $N$ different series in the conjugacy
class,  $\varphi_1,\ldots,\varphi_N$, we get
\[r_i(z)=z\varphi_i(z^{-1})=mz+a_1c_i^{N_1}z^{1-N_1/N}+a_2c_i^{N_2}z^{1-N_2/N}+a_3c_i^{N_3}z^{1-N_3/N}
+ \cdots\]
where $c_1,\ldots,c_N$ are the $N$ complex
roots of $x^N=1$.

\para

\noindent
Under these conditions, we introduce the following definition of branch.

\begin{definition}\label{D-infinitybranch}
 An {\em
infinity branch of an affine plane curve ${\cal C}$}  associated to the infinity
point $P=(1:m:0),\,m\in {\Bbb C}$, is  a set
$\displaystyle B=\bigcup_{j=1}^N L_j$, where $L_j=\{(z,r_j(z))\in
{\Bbb C}^2: \,z\in {\Bbb C},\,|z|>M\}$,\,  $M\in
{\Bbb R}^+$, and
\[r_j(z)=z\varphi_j(z^{-1})=mz+a_1c_j^{N_1}z^{1-N_1/N}+a_2c_j^{N_2}z^{1-N_2/N}+a_3c_j^{N_3}z^{1-N_3/N}
+ \cdots\] where $N, N_i\in {\Bbb N},\,\,i=1,\ldots$,  $0<N_1<N_2<\cdots$, and $c_j^N=1,\,\,j=1,\ldots,N$.
 The subsets $L_1,\ldots,L_N$ are
called the {\em leaves} of the infinity branch $B$.
 \end{definition}

 \para

\begin{remark} \label{R-conjugation} We observe that: \begin{enumerate}
\item An infinity branch  is uniquely determined  from one leaf, up to conjugation. That is, if
$\displaystyle B=\bigcup_{i=1}^NL_i$, where $L_i=\{(z,r_i(z))\in {\Bbb
C}^2: \,z\in {\Bbb C},\,|z|>M_i\}$, and
$$ r_i(z)=z\varphi_i(z^{-1})=mz+a_1z^{1-N_1/N}+a_2z^{1-N_2/N}+a_3z^{1-N_3/N} + \cdots$$
then $r_j=r_i,\,j=1,\ldots,N$, up to conjugation; i.e.
\[r_j(z)=z\varphi_j(z^{-1})=mz+a_1c_j^{N_1}z^{1-N_1/N}+a_2c_j^{N_2}z^{1-N_2/N}+a_3c_j^{N_3}z^{1-N_3/N}
+ \cdots\] where
$N, N_i\in\mathbb{N}$, and $c_j^N=1,\,\,j=1,\ldots,N$.
\item We may represent $L_i=\{(z,r_i(z))\in {\Bbb
C}^2: \,z\in {\Bbb C},\,|z|>M\}$,\,$i=1,\ldots,N$, where $M:=\max\{M_1,\ldots,M_N\}$.
\item  By abuse of notation, we say that $N$ is  {\em  the ramification index} of the branch $B$, and we  write
$\nu(B)=N$. Note that $B$ has $\nu(B)$  leaves.
\end{enumerate}
\end{remark}

\para
Let $\psi(t):=\varphi(t^N)$, where $\varphi(z)$ is a
series expansion for a solution  of $g(y,z)=0$. Observe that $(1:\psi(t):t^N)$ is a local projective parametrization, with
center at $P$, of the projective curve ${\cal C}^*$. Thus, from $\psi_i(t):=\varphi_i(t^N),\,\,i=1,\ldots,N$ ($\varphi_i$ are the $N$ different series in the conjugacy
class of $\varphi$),   we obtain
  $N$ equivalent local projective parametrizations,
$(1:\psi_i(t):t^N)$ (note that they are equivalent since   $\varphi_1,\ldots\varphi_N$ belong to the same
conjugacy class). Therefore, the leaves of $B$ are all
associated to a unique infinity place.

\para

Conversely, from a given  infinity place defined  by
a local projective pa\-ra\-me\-tri\-zation $(1:\psi(t):t^N)$ (see Theorem 2.5.3 in \cite{SWP}), we obtain $N$
Puiseux series, $\varphi_j(t)=\psi(c_j{t}^{1/N})$,
$c_j^N=1$, that provide different expressions
$r_j(z)=z\varphi_j(z^{-1}),\,j=1,\ldots,N$. Hence,  the infinity
branch $B$ is defined by the leaves
$L_i=\{(z,r_i(z))\in {\Bbb C}^2:
\,z\in {\Bbb C},\,|z|>M\},\,i=1,\ldots,N.$

\para

From the above discussion, we deduce that there exists a
one-to-one relation between infinity places and infinity branches.
In addition, we can say that each infinity branch is associated to
a unique infinity point given by the center of the corresponding infinity place. Reciprocally,  taking into account the above construction, we get that every infinity point has associated, at least, one
infinity branch. Hence, every algebraic plane curve has, at least, one infinity
branch. Furthermore,  every algebraic plane curve has a finite number of branches.

\para

In the following, we introduce the notions of convergent branches
and approaching curves. Intuitively speaking, two infinity branches
converge if they get closer  as they tend to infinity. This concept
will allow us to analyze whether two curves approach each other. For
further details see \cite{paper1}.

\begin{definition}\label{D-distance0}
Two
infinity branches, $B$ and $\overline{B}$, are convergent if there
exist two leaves   $L=\{(z,r(z))\in {\Bbb C}^2:\,z\in {\Bbb C},\,|z|>M\} \subset B$ and
$\overline{L}=\{(z,\overline{r}(z))\in {\Bbb C}^2:\,z\in {\Bbb
C},\,|z|>\overline{M}\}\subset \overline{B}$  such that  $\lim_{z\rightarrow\infty} (\overline{r}(z)-r(z))=0.$ In this case, we say that the leaves $L$ and $\overline{L}$ converge.
\end{definition}

Observe that   two convergent infinity branches are associated to
the same infinity point (see Remark 4.5 in \cite{paper1}).

\para

In the following lemma, we characterize the convergence of two given infinity branches (see Lemma 4.2, and Proposition 4.6 in \cite{paper1}).

\begin{lemma}\label{L-DistVertical} The following statements hold:
\begin{itemize}\item Two  leaves $L=\{(z,r(z))\in {\Bbb C}^2:\,z\in {\Bbb C},\,|z|>M\}$ and
$\overline{L}=\{(z,\overline{r}(z))\in {\Bbb C}^2:\,z\in {\Bbb
C},\,|z|>\overline{M}\}$ are convergent if and only if the terms
with non negative exponent in the series $r(z)$ and
$\overline{r}(z)$ are the same.
\item  Two infinity branches $B$ and $\overline{B}$ are convergent if and
only if for each leaf $L\subset B$ there exists a leaf $\overline{L}\subset
\overline{B}$ convergent with $L$, and reciprocally.\end{itemize}\end{lemma}

\para

Note that two convergent branches may be contained in the same curve
or they may belong to different curves. In this  second case, we will
say that {\it these curves approach each other}. More precisely, we have the following definition.

\begin{definition}\label{D-distance1}
Let ${\cal C} $ be an algebraic plane curve over ${\Bbb C}$ with an
infinity branch $B$. We say that a  curve ${\overline{{\cal C}}}$
{\em approaches} ${\cal C}$ at its infinity branch $B$ if there
exists one leaf $L=\{(z,r(z))\in {\Bbb C}^2:\,z\in {\Bbb
C},\,|z|>M\}\subset B$ such that
$\lim_{z\rightarrow\infty}d((z,r(z)),\overline{\cal C})=0.$
\end{definition}

\para

In the following, we state some important results concerning two curves that approach each other. These results are proved in \cite{paper1} (see Theorem 4.11 and Corollary 4.13).

\begin{theorem}\label{T-curvas-aprox}
Let ${\cal C}$
be a plane algebraic curve over $\Bbb C$ with an infinity branch
$B$. A plane algebraic curve ${\overline{{\cal C}}}$ approaches
${\cal C}$ at $B$ if and only if ${\overline{{\cal C}}}$ has an
infinity branch, $\overline{B}$, such that $B$ and $\overline{B}$
are convergent.
\end{theorem}

\para

\begin{remark}\label{R-approaching-curves}
\begin{enumerate}
\item It holds that ``{\em proximity}'' is a
symmetric relation; that is, ${\overline{{\cal C}}}$ approaches ${\cal C}$ at
some infinity branch $B$ if and only if ${\cal C}$ approaches
${\overline{{\cal C}}}$ at some infinity branch $\overline{B}$. In the following, we
say that ${\cal C}$ and ${\overline{{\cal C}}}$ approach each other
or that they are {\em approaching curves}.
\item Two approaching curves have
  a common infinity point.
\item  ${\overline{{\cal C}}}$ {\em
approaches} ${\cal C}$ at an infinity branch $B$ iff for every leaf $L=\{(z,r(z))\in {\Bbb C}^2:\,z\in {\Bbb
C},\,|z|>M\}\subset B$, it holds that
$\displaystyle\lim_{z\rightarrow\infty}d((z,r(z)),\overline{\cal
C})=0$.
\end{enumerate}
\end{remark}

\para

\begin{corollary}\label{C-approaching-curves}
Let $\cal C$ be an algebraic plane curve with an infinity branch
$B$. Let ${\overline{{\cal C}}}_1$ and ${\overline{{\cal C}}}_2$ be two
different curves that approach $\cal C$ at $B$. Then:
\begin{enumerate}
\item  ${\overline{{\cal C}}}_i$ has an infinity branch
$\overline{B_i}$ that converges with $B$, for $i=1,2$.
\item $\overline{B_1}$ and $\overline{B_2}$ are convergent. Then,   ${\overline{{\cal C}}}_1$
and ${\overline{{\cal C}}}_2$ approach each other.
\end{enumerate}\end{corollary}

\para

Taking into account that an infinity branch $B$  is uniquely
determined  from one leaf, up to conjugation (see statement 1 in
Remark \ref{R-conjugation}), and that the results stated above hold
for any leaf of $B$, for the sake of simplicity, in the following,
 $$B=\{(z,r(z))\in {\Bbb C}^2:\,z\in {\Bbb
C},\,|z|>M\}$$  stands for the infinity branch whose leaves are
obtained by conjugation on
$$r(z)=mz+a_1z^{1-N_1/N}+a_2z^{1-N_2/N}+a_3z^{1-N_3/N} + \cdots,\quad a_i\not=0,\, \forall i\in {\Bbb N} $$
$N,N_i\in {\Bbb N},\,\,i=1,\ldots$, and $0<N_1<N_2<\cdots$.  We also will prove that the results obtained throughout the paper hold for any leaf.

\section{Asymptotes and perfect curves}\label{S-asymptotes}

Given an algebraic plane curve $\cal C$ and a infinity branch $B$ of
$\cal C$, in \cite{paper1} we analyze whether $\cal C$ can  be
approached  at $B$ by a new curve ${\overline{\cal C}}$. Intuitively speaking, if $\cal C$ is
approached  at $B$ by  ${\overline{\cal C}}$,  and $\deg({\overline{\cal
C}})<\deg({\cal C})$,  one may say that $\cal
C$ {\it degenerates}, since $\cal C$ behaves at the infinity as a curve of
less degree.

\para

For instance, a hyperbola is a curve of degree 2 that has two real
asymptotes, which implies that the hyperbola  degenerates, at
the infinity, in two lines. The behavior of an ellipse is
similar; in this case, the infinity branches are complex but they can
also be approached by (complex)   lines. However, the
asymptotic behavior of a parabola is   different, since at the infinity, the parabola
cannot be approached by any   line. This  motivates the
following definition.

\begin{definition}\label{D-perfect-curve}
A curve of degree $d$ is a {\em perfect curve} if it cannot be approached by
any curve of degree less than $d$.
\end{definition}

A curve that is not perfect can be approached by other curves of
less degree. If these curves are perfect, we call them
{\it g-asymptotes}. More precisely, we have the following definition.

\begin{definition}\label{D-asymptote}
Let ${\cal C}$ be a curve with an infinity branch $B$. A {\em
g-asymptote} (generalized asymptote) of ${\cal C}$ at $B$ is a
perfect curve that approaches ${\cal C}$ at $B$.
\end{definition}

Note that the notion of {\em
g-asymptote} is similar to the classical concept of asymptote. The
difference is that a g-asymptote does not have to be a line, but a
perfect curve. Actually, it is a generalization, since every line is
a perfect curve (this fact follows from Definition
\ref{D-perfect-curve}). Throughout the paper we  refer to
{\em g-asymptote} simply as {\it asymptote}.

\para

\begin{figure}[h]
$$
\begin{array}{cc}
\psfig{figure=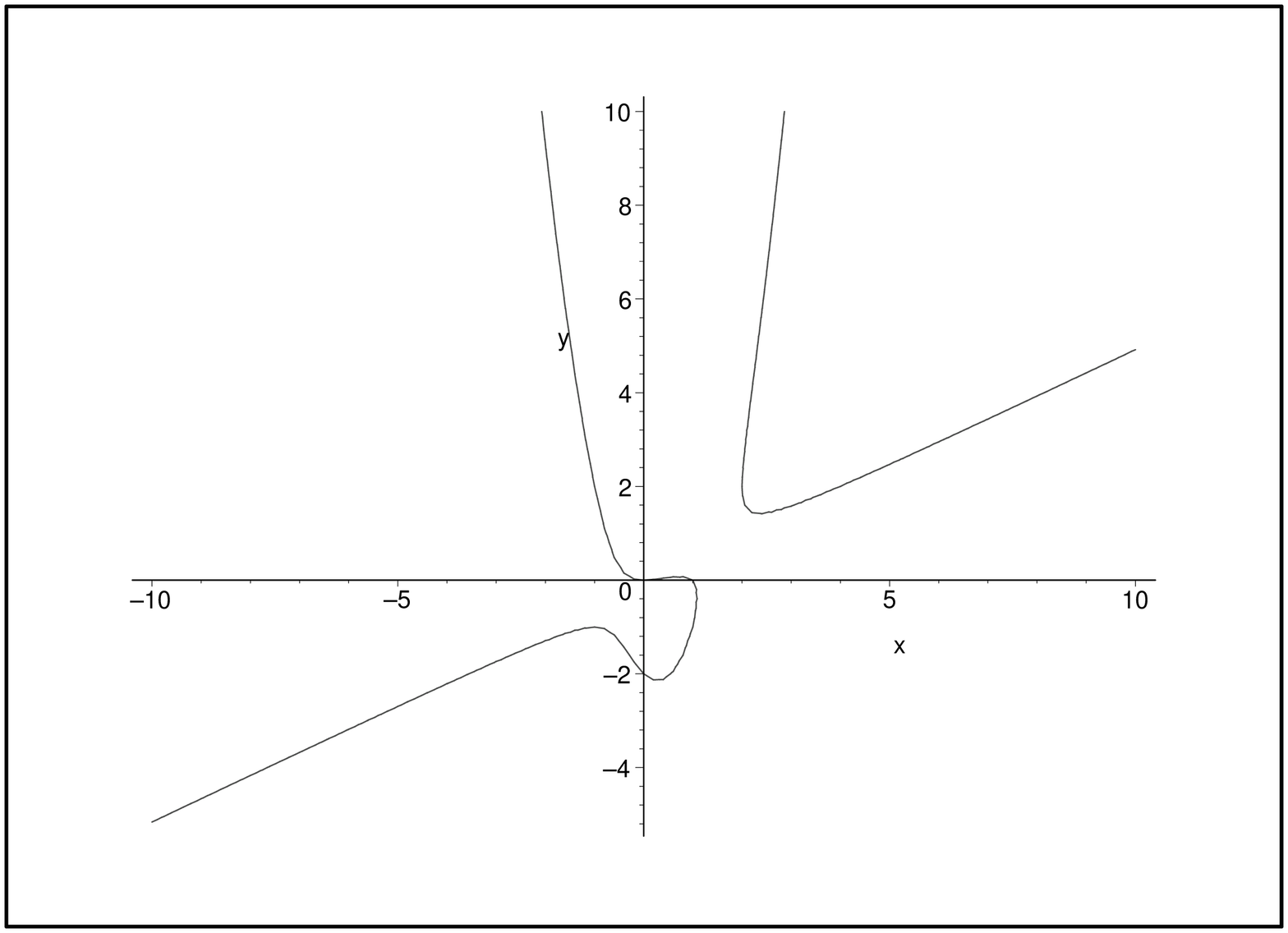,width=4.3cm,height=4.3cm,angle=270} &
\psfig{figure=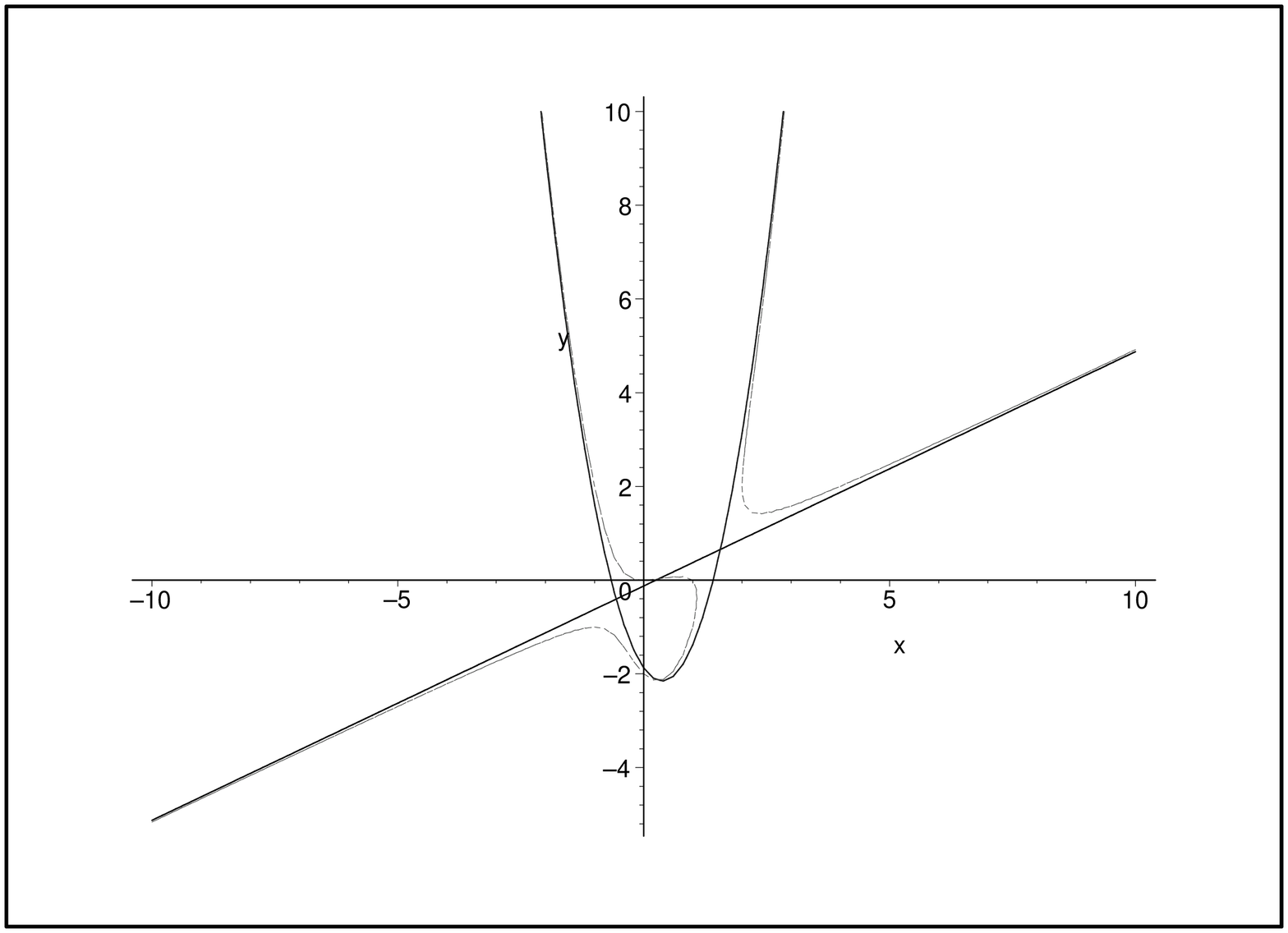,width=4.3cm,height=4.3cm,angle=270}
\end{array}
$$ \caption{Curve $\cal C$ (left)  approached by a parabola and a line (right).}\label{F-ejemplo-asintotas}
\end{figure}

In order to clarify this notion, let us consider a plane curve  $\cal C$   defined by the irreducible polynomial
$$f(x,y)=-yx-y^2-x^3+2x^2y+x^2-2y\in {\Bbb R}[x,y].$$
$\cal C$ has degree 3, and two infinity branches. In Figure
\ref{F-ejemplo-asintotas}, one can check that these infinity
branches are approached  by the parabola
$y-2x^2+3/2x+15/8=0$, and the line
$y-x/2+1/8=0$.
\vspace*{-0.5cm}

Later we see  that other plane curves of degree 3,
like $y-x^3=0$ or $y^2-x^3=0$, cannot be approached by any curve of
degree less than 3. That is, they are perfect curves.

\begin{remark}\label{R-minimal-degree}
The degree of an asymptote is less or equal than the
degree of the curve it approaches. In fact, an asymptote of a curve
$\cal C$ at a branch $B$ has minimal degree among all the curves
that approach $\cal C$ at $B$. Indeed, let ${\cal D}$
be an asymptote of $\cal C$ at $B$ and let ${\cal D}'$ be another
curve that approaches $\cal C$ at $B$. From Corollary
\ref{C-approaching-curves}, ${\cal D}'$ approaches ${\cal D}$, and since  ${\cal D}$ is
  perfect, we conclude that $\deg({\cal D}')\geq \deg({\cal D})$.
\end{remark}

\para

 In the following, we show that every infinity branch
of a given algebraic plane curve has, at least, one asymptote (see Theorem
\ref{T-constr-asintota}). In order to prove this property, we first
need to show some previous results. For this purpose, let   ${\cal C}$ be  a plane curve and $B=\{(z,r(z))\in {\Bbb C}^2:\,z\in {\Bbb C},\,|z|>M\}$ an   infinity
branch of  ${\cal C}$  associated to
$P=(1:m:0)$. From Definition
\ref{D-infinitybranch}, we have that
 \begin{equation}\label{Eq-inf-branchN}r(z)=mz+a_1z^{-N_1/N+1}+\cdots
+a_kz^{-N_k/N+1}+a_{k+1}z^{-N_{k+1}/N+1}+\cdots\end{equation}  where
$a_1,a_2,\ldots\in\mathbb{C}\setminus \{0\}$,\,$m\in {\Bbb C}$,
$N, N_1,N_2\ldots\in\mathbb{N}$, and $0<N_1<N_2<\cdots$. In addition,
let $N_k\leq N <N_{k+1}$, i.e. the terms $a_jz^{-N_j/N+1}$
with $j\geq k+1$ have negative exponent. Note that   $\nu(B)=N$.

\para

\begin{lemma}\label{L-case3}
Let $\cal C$ be a plane curve over $\Bbb C$ containing an infinity
branch   of the form in
(\ref{Eq-inf-branchN}), and let $f(x,y)\in {\Bbb R}[x,y]$ be the
irreducible polynomial defining implicitly   $\cal C$. It holds that
$(y-mx)^N$ divides $f_d(x,y)$.
\end{lemma}

\noindent\textbf{Proof:}  Let $r_1,\ldots,r_N$ be the conjugates of
$r$. That is,
\[r_j(z)=z\varphi_j(z^{-1})=mz+a_1c_j^{N_1}z^{1-N_1/N}+a_2c_j^{N_2}z^{1-N_2/N}+a_3c_j^{N_3}z^{1-N_3/N}
+ \cdots,\] where $c_j^N=1,\,\,j=1,\ldots,N$.  Now, we consider
$$h(x,y):=\prod_{i=1}^N(y-r_i(x)).$$
Note that the terms with maximum exponent in $h(x,y)$ are given by $(y-mx)^N$.  We denote $h_1(x,y):=(y-mx)^N$.
% (we recall  that these series are associated to the different leaves of the branch $B$; see Definition \ref{D-infinitybranch})

\para
\noindent
 In addition, if we see $f(x,y)$ and
$h(x,y)$ as  polynomials in the variable $y$, and we denote it by
$f_x(y),h_x(y)\in\mathbb{C}\ll x\gg[y]$, it holds that $h_x$ divides
$f_x$ (w.r.t the variable $y$). Indeed: let $\varphi_1,\ldots,\varphi_N$ be the
series expansions for solutions of $g(y,z)=0$, where
$g(y,z)=F(1:y:z)$. Now, we see $g(y,z)$ as a polynomial in the
variable $y$, and we denote it by $g_z(y)\in\mathbb{C}\ll z\gg[y]$.
From Theorem 4.2 in \cite{verger}, we deduce that
$(y-\varphi_i(z))$ divides $g_z(y)$, for  $i=1,\ldots,N$. Hence,
$x(y/x-\varphi_i(z/x))$ divides $x^dF(1:y/x:z/x)=F(x:y:z)$  w.r.t the variable $y$.
Now, we set $z=1$,  and we deduce that
$(y-x\varphi_i(x^{-1}))=(y-r_i(x))$ divides $F(x:y:1)=f_x(y)$. Furthermore, since the factors of $h(x,y)$ are all different (they are obtained by conjugation), we get  that
$$f_x(y)=h_x(y)p_x(y),\quad \mbox{where}\quad p_x(y)\in\mathbb{C}\ll x\gg[y].$$

\noindent Now, we prove that  $(y-mx)^N$ divide $f_d(x,y)$. For this
purpose, we write
\[h_x(y)=h_1(x,y)+\overline{q}_1(x,y),\quad
p_x(y)=p_1(x,y)+\overline{q}_2(x,y),\quad \overline{q}_i\in\mathbb{C}\ll x\gg[y],\] where
$h_1(x,y)=(y-mx)^N$, and $p_1(x,y):=\prod_{i=1}^r(y-m_ix)\,\in
{\Bbb C}[x,y],\,m_i\in {\Bbb C}$ are homogeneous polynomials. Note that
$\deg(h_1)=\deg_y(h_x)=N$,   $\deg(p_1)=\deg_y(p_x)=r$, and $\deg_y(\overline{q}_1)<N$,
$\deg_y(\overline{q}_2)<r$.  Then, we may write
\[z^N h_x(y/z)=(y-mxz)^N+zq_1(x,y,z),\quad z^r p_x(y/z)=\prod_{i=1}^r(y-m_ixz)+zq_2(x,y,z),\]
where $q_i\in\mathbb{C}\ll x\gg[y,z],\,i=1,2$. In addition, since
\[f(x,y)=f_d(x,y)+f_{d-1}(x,y)+f_{d-2}(x,y)+\cdots+f_0 \in {\Bbb R}[x,y],\]
where  $f_k,\,k=0,\ldots,d,$ are
homogeneous polynomials of degree $k$, we have that
\[z^d f_x(y/z)=z^df(x,y/z)=z^df_d(x,y/z)+z^df_{d-1}(x,y/z)+\cdots+z^d f_0=\]\[z^df_d(x,y/z)+zq(x,y,z),\]
where   $q\in {\Bbb R}[x,y,z]$. Observe that $f_d(x,y)=\prod_{i=1}^d(y-s_ix),\,s_i\in {\Bbb C}$, and then
$z^df_d(x,y/z)=\prod_{i=1}^d(y-s_ixz)$.\\

\noindent
Under these conditions,  since $f_x=h_xp_x$, and $d=N+r$, we have that
\[z^d f_x(y/z)=z^N h_x(y/z) z^r p_x(y/z) \]
which implies that
\[\prod_{i=1}^d(y-s_ixz)+zq(x,y,z)=((y-mxz)^N+zq_1(x,y,z))(\prod_{i=1}^r(y-m_ixz)+zq_2(x,y,z)).\]
That is, \[\prod_{i=1}^d(y-s_ixz)-(y-mxz)^N
\prod_{i=1}^r(y-m_ixz)=\]\[z(-q(x,y,z)+q_1(x,y,z)\prod_{i=1}^r(y-m_ixz)+q_2(x,y,z)(y-mxz)^N+zq_1(x,y,z)q_2(x,y,z)).\]
Since $z$ does not divide the left hand side of  the above equality,
we get that $$\prod_{i=1}^d(y-s_ixz)=(y-mxz)^N
\prod_{i=1}^r(y-m_ixz).$$ Thus, $(y-mxz)^N$ divides
$z^df_d(x,y/z)=\prod_{i=1}^d(y-s_ixz)$, and then $(y-mx)^N$ divides
$f_d(x,y)$.\hfill $\Box$
%Now, let $g_1(x,y)$ and $h_1(x,y)$
%be the terms with maximum exponent in $g(x,y)$ and $h(x,y)$,
%respectively (here we mean the global exponent, obtained by adding
%the exponents of the variables $x$ and $y$). Clearly, the form of
%maximum degree in $f(x,y)$ is given by $g_1(x,y)h_1(x,y)$, but we
%have seen above that $h_1(x,y)=(y-mx)^n$, so the result
%follows.

\para

\para

\noindent
In the following, we write equation (\ref{Eq-inf-branchN}) defining a branch $B$ as
 \begin{equation}\label{Eq-inf-branchn}r(z)=mz+a_1z^{-n_1/n+1}+\cdots
+a_kz^{-n_k/n+1}+a_{k+1}z^{-N_{k+1}/N+1}+\cdots\end{equation} where
$\gcd(N,N_1,\ldots,N_k)=b$,\,$N_j=n_jb,\,N=nb,\,\,j=1,\ldots,k$. That is, we have simplified the non negative
exponents such that $\gcd(n,n_1,\ldots,n_k)=1$. Note that $0<n_1<n_2<\cdots$,  $n_k\leq n$, and $N<n_{k+1}$, i.e.  the terms $a_jz^{-N_j/N+1}$
with $j\geq k+1$ have negative exponent. We denote these terms as
\[ A(z):=\sum_{\ell=k+1}^\infty
a_{\ell}z^{-q_{\ell}},\,\quad
q_{\ell}=-N_{\ell}/N+1\in\mathbb{Q}^+,\,\,\,\ell \geq k+1.\]

Under these conditions, we introduce the definition of degree of a branch $B$ as follows:

\para

\begin{definition}\label{D-degreebranch} Let $B=\{(z,r(z))\in
{\Bbb C}^2: \,z\in {\Bbb C},\,|z|>M\}$ defined by  (\ref{Eq-inf-branchn})  an
infinity branch associated to   $P=(1:m:0),\,m\in {\Bbb C}$. We say that $n$ is de degree of $B$, and we denote it by $\deg(B)$.
\end{definition}

\begin{proposition}\label{C-case3}
Let $\overline{{\cal C}}$ be a curve that approaches ${\cal C}$ at
its infinity branch $B$. Let $\overline{f}\in\mathbb{R}[x,y]$ be the
implicit polynomial of  $\overline{{\cal C}}$ . Then, $(y-mx)^n$
divides the homogeneous form of maximum degree of
$\overline{f}(x,y)$.
\end{proposition}

\noindent\textbf{Proof:} Using Theorem \ref{T-curvas-aprox}, we get
that  $\overline{{\cal C}}$ has an infinity branch
$\overline{B}=\{(z,\overline{r}(z))\in {\Bbb C}^2:\,z\in {\Bbb
C},\,|z|>\overline{M}\}$ convergent with $B=\{(z,r(z))\in {\Bbb
C}^2:\,z\in {\Bbb C},\,|z|>M\}$. From Lemma \ref{L-DistVertical}, we
deduce that the terms with non negative exponent in the series
$r(z)$ and $\overline{r}(z)$ are the same, and hence $\overline{B}$
is a branch of degree $n$ of the form given in
(\ref{Eq-inf-branchn}).  Thus, Lemma \ref{L-case3} states
that the homogeneous form of maximum degree of $\overline{f}(x,y)$
is divided by $(y-mx)^{\nu(\overline{B})}$. Now, the result follows
taking into account that, always, $\deg(\overline{B})=n\leq
\nu(\overline{B})$.\hfill $\Box$

\para

\begin{remark}\label{R-grado-caso3}
From Lemma \ref{L-case3}, we deduce that a curve $\cal C$ containing
an infinity branch $B$ of degree $n$  has degree at least $\nu(\overline{B})\geq n$.
Furthermore,  from Proposition \ref{C-case3}, we also have that $\deg(\overline{\cal C})\geq n$ for any curve $\overline{\cal C}$ approaching $\cal C$ at $B$.
\end{remark}

\subsection{Construction of an asymptote}

Taking into account the results presented  above, we have  that any
curve $\overline{{\cal C}}$ approaching ${\cal C}$ at $B$ should
have an infinity branch
 $\overline{B}=\{(z,\overline{r}(z))\in {\Bbb
C}^2:\,z\in {\Bbb C},\,|z|>\overline{M}\}$ such that the terms with
non negative exponent in $r(z)$ and $\overline{r}(z)$ are the same.
In the simplest case, if $A=0$ (i.e. there are not terms with negative
exponent; see equation (\ref{Eq-inf-branchn})), we obtain
\begin{equation}\label{Eq-inf-branch3}
\tilde{r}(z)=mz+a_1z^{-n_1/n+1}+a_2z^{-n_2/n+1}+\cdots
+a_kz^{-n_k/n+1},
\end{equation}
where $a_1,a_2,\ldots\in\mathbb{C}\setminus \{0\}$,\,$m\in {\Bbb
C}$, $n,n_1,n_2\ldots\in\mathbb{N}$, $\gcd(n,n_1,\ldots,n_k)=1$, and $0<n_1<n_2<\cdots$.
Note that $\tilde{r}$ has the same terms with non negative exponent that
$r$, and $\tilde{r}$ does not have terms with negative exponent.

\para

Let $\widetilde{{\cal C}}$ be the plane curve containing the branch
$\widetilde{B}=\{(z,\tilde{r}(z))\in {\Bbb C}^2:\,z\in {\Bbb
C},\,|z|>\widetilde{M}\}$ (note that $\widetilde{{\cal C}}$  is
unique since two different algebraic curves have finitely many
common points). Observe that
\begin{equation}\label{Eq-parametric-case1}\widetilde{{\cal Q}}(t)=(t^n,mt^n+a_1t^{n-n_1}
+\cdots +a_kt^{n-n_k})\in {\Bbb C}[t]^2,\end{equation}
where  $n,n_1,\ldots,n_k\in\mathbb{N}$, $\gcd(n,n_1,\ldots,n_k)=1$, and $0<n_1<\cdots <n_k$,
is a polynomial  parametrization of $\widetilde{{\cal C}}$, and it is proper (see Lemma \ref{L-proper-param}). In Theorem \ref{T-constr-asintota}, we  prove that
$\widetilde{{\cal C}}$ is an asymptote of ${\cal C}$ at $B$.

%For this purpose, we first prove  that $\widetilde{{\cal Q}}$ is a proper  (i.e. invertible) parametrization.

\para

\begin{lemma}\label{L-proper-param} The parametrization given in (\ref{Eq-parametric-case1}) is proper
 (i.e. invertible).
\end{lemma}

\noindent\textbf{Proof:} Let us assume that  $\widetilde{{\cal Q}}$ is not proper.
Then,
 there exists $R(t)\in {\Bbb C}[t]$, with $\deg(R)=r>1$,
and ${\cal Q}(t)=(q_1(t),q_2(t))\in {\Bbb C}[t]^2$, such that ${\cal
Q}(R)=\widetilde{{\cal Q}}$ (see  \cite{Binder}, \cite{Hong} or \cite{Manocha}).
In particular, we have that $q_1(R(t))=t^n$, which implies that
$$q_1(t)=(t-R(0))^k,\quad \mbox{and}\quad
R(t)=t^r+R(0),\quad rk=n.$$ Let us consider
$R^{\star}(t)=R(t)-R(0)=t^r\in {\Bbb C}[t]$, and
\[{\cal Q}^{\star}(t)={\cal Q}(t+R(0))=(t^k,{q^{\star}_2}(t))=(t^k, c_0+c_1t+c_2t^2+\ldots+c_ut^u) \in {\Bbb C}[t]^2.\]
Then, ${\cal Q}^{\star}(R^{\star})={\cal Q}(R)=\widetilde{{\cal Q}}$, and in
particular   ${q^{\star}_2}(R^\star)={q^{\star}_2}(t^r)=mt^n+a_1t^{n-n_1}+a_2t^{n-n_2}+\cdots +a_kt^{n-n_k}$.
That is,
  $$c_0+c_1t^r+c_2t^{2r}+\ldots+c_ut^{ur}=mt^n+a_1t^{n-n_1}+a_2t^{n-n_2}+\cdots +a_kt^{n-n_k}.$$
  From this equality, and taking into account that $r$ divides $n$ (recall that $rk=n$), we deduce that   $r$ divides $n_j,\,j=1,\ldots,k$.   This is impossible, because $r>1$, and $\gcd(n,n_1,\ldots,n_k)=1$. Therefore, we conclude that  $\widetilde{{\cal Q}}$ is proper.\hfill
$\Box$

\para

The next result states a property concerning the implicit polynomial
of $\widetilde{{\cal C}}$ (compare with Lemma \ref{L-case3}).

\para

\begin{lemma}\label{L-case1}
Let $\widetilde{{\cal C}}$ be the plane curve containing the infinity
branch given in (\ref{Eq-inf-branch3}). Let $\tilde{f}(x,y)\in
{\Bbb R}[x,y]$ be the implicit polynomial defining  $\widetilde{{\cal C}}$. It holds that the homogeneous form
of maximum degree of $\tilde{f}(x,y)$ is $(y-mx)^n$.
\end{lemma}
\noindent\textbf{Proof:} First, we consider the polynomial proper parametrization (see Lemma
\ref{L-proper-param}) defining $\widetilde{{\cal C}}$, and introduced in  (\ref{Eq-parametric-case1}):
\[\widetilde{{\cal Q}}(t)=(t^n,mt^n+a_1t^{n-n_1}
+\cdots +a_kt^{n-n_k})\in {\Bbb C}[t]^2.\]Now, we distinguish two different cases:

\begin{itemize}
\item[1.] If $m=0$, i.e. $B$ is associated to the
infinity point $P=(1:0:0)$, we apply the results in  \cite{SWP}
(see Chapter 4), and one has that, up to constants in ${\Bbb R}\setminus\{0\}$,
\[%\label{Eq-implicit-case1}
\tilde{f}(x,y)=\resultant_t(x-t^n,y-p(t))=\prod_{i=1}^n(y-p(\alpha_i))
\]
where
$p(t)=mt^n+a_1t^{n-n_1}+\cdots +a_kt^{n-n_k}$, and $\alpha_1,\ldots,\alpha_n$ are the $n$ roots of the equation
$x-t^n=0$. Hence, since $\deg(p)=n-n_1$, we get that  the
maximum exponent of $p(\alpha_i)$ is $(n-n_1)/n<1$  and then, the  form of
maximum degree of $\tilde{f}(x,y)$ is $y^n$.
\item[2.] Let  $m\not=0$, and then $B$ is associated to the
infinity point  $P=(1:m:0)$. In this case,
we  apply the linear change of variables, $x=X-mY,\, y=mX+Y$, and  the infinity point moves to $(1:0:0)$. By applying   case 1, we get that  the
homogeneous form of maximum degree of  $\tilde{f}(X-mY,mX+Y)$
is $Y^n$. Finally, undoing the change, we get that the homogeneous form
of maximum degree of $\tilde{f}(x,y)$ is $(y-mx)^n$.\hfill $\Box$\end{itemize}

\begin{remark}\label{R-grado-caso1}
From Lemma \ref{L-case1}, we deduce that $\deg(\widetilde{{\cal C}})=n$.
\end{remark}

\begin{theorem}\label{T-constr-asintota}
The curve $\widetilde{{\cal C}}$ is an asymptote of $\mathcal{C}$ at $B$.
\end{theorem}

\noindent\textbf{Proof:} Taking into account the construction of
$\widetilde{{\cal C}}$, we have that $\widetilde{{\cal C}}$
approaches $\mathcal{C}$ at $B$. Therefore, we only  need to show
that $\widetilde{{\cal C}}$ is perfect, i.e.  that $\widetilde{{\cal
C}}$ cannot be approached by any curve with degree less than
$\deg(\widetilde{{\cal C}})$.

\para

For this purpose, we first note that $\widetilde{{\cal C}}$ admits the
polynomial parametrization given by the form in  (\ref{Eq-parametric-case1}). Then, using the results in \cite{Manocha}, we deduce that the unique infinity branch of $\widetilde{{\cal C}}$ is $\widetilde{B}$. In addition, we observe that by construction, $\widetilde{B}$ and $B$ are convergent.
%In \cite{Manocha}, it is proved that a curve is polynomial if and only
%if it has only  one infinity point $P$, and the corresponding projective
%curve has   only  one place centered at $P$. That is, a
%curve is polynomial if and only if it has  only one infinity branch. Hence, we deduce that the unique infinity branch of $\widetilde{{\cal C}}$ is
%$\widetilde{B}$, and we note that by construction, $\widetilde{B}$ and
%$B$ are convergent.

\para

Under these conditions, let us consider  a plane curve, $\overline{\cal C}$, that approaches $\widetilde{{\cal C}}$
at $\widetilde{B}$. From Theorem \ref{T-curvas-aprox}, we get that $\overline{\cal C}$
  has an infinity branch $\overline{B}$ convergent with $\widetilde{B}$. Since   $\widetilde{B}$ and
$B$ are convergent, from Corollary \ref{C-approaching-curves}, we deduce that $\overline{B}$ and $B$ are  convergent which implies  that $\overline{\cal C}$ approaches $\cal C$ at $B$. Now, from
Remarks \ref{R-grado-caso3} and \ref{R-grado-caso1}, we deduce that
$\deg(\overline{\cal C})\geq n=\deg(\widetilde{{\cal C}})$. Therefore, we conclude that $\widetilde{{\cal C}}$ is perfect.\hfill $\Box$

\para

We have shown that, for any infinity branch $B$ of a plane
curve $\cal C$, there always exists an asymptote that approaches
$\cal C$ at $B$. Furthermore, we have provided a method to obtain it. From these results, we obtain the following algorithm that  computes
an asymptote for each infinity branch of a given plane curve.

 \para
 We assume that we have  prepared the input curve $\cal C$, such that  by means  of a suitable linear change of coordinates,  $(0:1:0)$ is not an infinity point of $\cal C$.

\para

\begin{center}
\fbox{\hspace*{2 mm}\parbox{13.2cm}{ \vspace*{2 mm} {\bf Algorithm
{\sf Asymptotes Construction.}} \vspace*{0.2cm}

\noindent Given an implicit algebraic plane curve $\cal C$ over $\Bbb C$, the algorithm computes one
asymptote for each of its infinity branches.\vspace*{0.1cm}

\noindent
1. Compute the infinity points of $\cal C$. Let $P_1,...,P_n$ be these points.\vspace*{0.1cm}
%\item[2.] If $P_j=(0:1:0)$ for some $j\in \{1,\ldots,n\}$, apply a linear change of variables such that $P_j$ is moved to $(1:m:0)$, for some $m\in {\Bbb C}$, and $(0:1:0)$ is not an infinity point of the new curve.

\noindent
2. For each $P_i:=(1:m_i:0)$ do:
\begin{itemize}
\item[2.1.] Compute the infinity branches of  $\cal C$ associated
to $P_i$. Let $B_j=\{(z,r_j(z))\in {\Bbb C}^2:\,z\in {\Bbb C},\,|z|>M_j\}$, $j=1,\ldots, s_i,$  be these branches, where $r_j$ is written as in equation (\ref{Eq-inf-branchn}). That is,
\[r_j(z)=m_iz+a_{1,j}z^{-n_{1,j}/n_j+1}+\cdots
+a_{k_j,j}z^{-n_{k_j,j}/n_j+1}+A_j(z),\] \[ A_j(z)=\sum_{\ell=k_j+1}^\infty
a_{\ell, j}z^{-q_{\ell, j}},\,\quad q_{\ell,j}=-N_{\ell,j}/N_j+1\in\mathbb{Q}^+,\,\,\,\ell \geq k_j+1,\]
$a_{1,j},a_{2,j},\ldots\in\mathbb{C}\setminus \{0\}$,\,
$n_j,n_{1,j},\ldots\in\mathbb{N}$,  $0<n_{1,j}<n_{2,j}<\cdots$, $n_{k_j}\leq n_j$, $N_j<n_{k_j+1}$, and $\gcd(n_j,n_{1,j},\ldots,n_{k_j,j})=1$.
\item[2.2.] For each branch $B_j,\,j=1,\ldots, s_i$ do:
\begin{itemize}
\item[2.2.1.] Consider $\tilde{r}_j$ as in equation (\ref{Eq-inf-branch3}). That is, \[\tilde{r}_j(z)=m_iz+a_{1,j}z^{-n_{1,j}/n_j+1}+\cdots
+a_{{k_j},j}z^{-n_{{k_j},j}/n_j+1}\]
Note that $\tilde{r}$ has the same terms with non negative
exponent that $r$, and $\tilde{r}$ does not have terms with negative exponent.
%$\tilde{r}_j(z)$ having the same terms with non negative
%exponent that $r_j(z)$, and  not having terms with negative exponent.
%by canceling the terms with negative exponent in $r_j(z)$.
\item[2.2.2.] {\sf Return} the asymptote $\widetilde{\cal C}_j$ defined by the proper  parametrization (see Lemma \ref{L-proper-param}),
$\widetilde{Q}_j(t)=(t^{n_j},\,\tilde{r}_j(t^{n_j}))\in {\Bbb C}[t]^2$,
%where $n_j$ is the
%common denominator of the exponents in $\tilde{r}_j(z)$,
and the implicit polynomial  (see Chapter 4 in \cite{SWP}):
 $$\tilde{f}_j(x,y)=\Resultant_t(x-t^{n_j},y-\tilde{r}_j(t^{n_j}))\in {\Bbb C}[x,y].$$
\end{itemize}
\end{itemize}
}\hspace{2 mm}}
\end{center}

\para

\noindent
{\it Correctness.}  The algorithm
{\sf Asymptotes Construction} outputs an asymptote $\widetilde{\cal C}$ that is independent of leaf chosen to define  the branch $B=\{(z,r(z))\in {\Bbb C}^2:\,z\in {\Bbb C},\,|z|>M\}$ (see Section 2). Indeed: let  $\widetilde{\cal C}$ be an asymptote obtained by the algorithm, and defined by the   proper parametrization  $\widetilde{Q}(t)=(t^{n},\,\tilde{r}(t^{n}))$, where
\[\tilde{r}(z)=mz+a_{1}z^{-n_{1}/n+1}+\cdots
+a_{k}z^{-n_{k}/n+1},\qquad \mbox{and}\]
\[r(z)=mz+a_{1}z^{-n_{1}/n+1}+\cdots
+a_{k}z^{-n_{k}/n+1}+A(z),\quad A=\sum_{\ell=k+1}^\infty
a_{\ell}z^{-q_{\ell}},\,\,
q_{\ell}\in\mathbb{Q}^+,\]%q_{\ell}=1-\frac{n_{\ell}}{N}
$a_{1},a_{2},\ldots\in\mathbb{C}\setminus \{0\}$,\,
$n,n_{1},\ldots\in\mathbb{N}$,  $0<n_{1}<n_{2}<\cdots$, and
$\gcd(n,n_{1},\ldots,n_{k})=1$,\,$N_j=n_jb,\,N=nb,\,\,j=1,\ldots,k$  (see equations (\ref{Eq-inf-branchn}) and (\ref{Eq-inf-branch3})). Now, let
\[r_{s}(z)=mz+a_{1}c_s^{n_{1}b}z^{-n_{1}/n+1}+\cdots
+a_{k}c_s^{n_{k}b}z^{-n_{k}/n+1}+A_s(z),\,\,
A_s=\sum_{\ell=k+1}^\infty
a_{\ell}c_s^{n_{\ell}}z^{-q_{\ell}},\]
where $c_s^{N}=1,\,\,s=1,\ldots,N$. That is,  $r_{s}=r$, up to
conjugation. Then, the parametrization obtained by algorithm using $r_s$ is
$\widetilde{Q}_s(t)=(t^{n},\,\tilde{r}_s(t^{n})),$
where
\[\tilde{r}_s(z)=mz+a_{1}c_s^{n_{1}b}z^{-n_{1}/n+1}+\cdots
+a_{k}c_s^{n_{k}b}z^{-n_{k}/n+1}.\]
Since
$$\widetilde{Q}(t)=(t^{n},\,\tilde{r}(t^{n}))=(t^{n}, mt^n+a_{1}t^{-n_{1}+n}+\cdots
+a_{k}z^{-n_{k}+n}),\quad \mbox{and}$$
$$\widetilde{Q}_s(t)=(t^{n},\,\tilde{r}_s(t^{n}))=(t^n, mt^n+a_{1}c_s^{n_{1}b}t^{-n_{1}+n}+\cdots
+a_{k}c_s^{n_{k}b}t^{-n_{k}+n}),$$ and taking into account that
$c_s^{N}=c_s^{nb}=1$, we deduce that $\widetilde{Q}_s(c_s^{b}t)=\widetilde{Q}(t).$
Therefore, both parametrizations, $\widetilde{Q}_s$ and
$\widetilde{Q}$, define the same asymptote $\widetilde{\cal C}$.

\para

In the following, we illustrate algorithm {\sf Asymptotes
Construction} with an example.

\begin{example}\label{Ej-asint-cuartica} Let ${\cal C}$ be the curve of degree $d=4$  defined  by the irreducible polynomial
$$f(x,y)=2y^3x-y^4+2y^2x-y^3-2x^3+x^2y+3\in {\Bbb R}[x,y].$$
We apply  algorithm {\sf Asymptotes Construction} to compute
the asymptotes of $\cal C$.

\begin{itemize}
\item[] \mbox{\sf {Step 1}:} We have that $f_4(x,y)=2y^3x-y^4$. Hence,  the infinity points are
$P_1=(1:2:0)$ and $P_2=(1:0:0)$.\\

We start by analyzing the point $P_1$:

\item[]  \mbox{\sf {Step 2.1}:} The only infinity branch associated to $P_1$
is $B_1=\{(z,r_1(z))\in {\Bbb C}^2:\,z\in {\Bbb C},\,|z|>M_1\}$, where
$$r_1(z)=2z+\frac{3z^{-3}}{8}-\frac{9z^{-4}}{64}+\frac{27z^{-5}}{512}-\frac{81z^{-6}}{4096}+\cdots $$
(we compute $r_1$ using the {\sf algcurves} package included in  the computer algebra system {\sf Maple}).
\item[]  \mbox{\sf {Step 2.2.1}:} We compute $\tilde{r}_1(z)$, and we have that
$\tilde{r}_1(z)=2z.$

\item[]  \mbox{\sf {Step 2.2.2}:} The parametrization of the asymptote $\widetilde{\cal C}_1$ is given by
$\widetilde{Q}_1(t)=(t,2t)\in {\Bbb R}[t]^2,$
and the polynomial defining implicitly $\widetilde{\cal C}_1$ is
$$\tilde{f}_1(x,y)=y-2x\in {\Bbb R}[x,y].$$

Now, we focus on the point $P_2$:
\item[]  \mbox{\sf {Step 2.1}:} The only infinity branch associated to $P_2$
is $B_2=\{(z,r_2(z))\in {\Bbb C}^2:\,z\in {\Bbb C},\,|z|>M_2\}$, where
$$r_2(z)=z^{2/3}-\frac{1}{3}+\frac{z^{-2/3}}{9}-\frac{2z^{-4/3}}{81}+\cdots.$$
\item[]  \mbox{\sf {Step 2.2.1}:} We obtain that
$\tilde{r}_2(z)=z^{2/3}-\frac{1}{3}.$

\item[]  \mbox{\sf {Step 2.2.2}:}  The parametrization of the asymptote $\widetilde{\cal C}_2$ is given by
$\widetilde{Q}_2(t)=(t^3,t^2-1/3)\in {\Bbb R}[t]^2,$
and the polynomial defining implicitly $\widetilde{\cal C}_2$ is
$$\tilde{f}_2(x,y)=-x^2+y^3+y^2+1/3y+1/27\in {\Bbb R}[x,y].$$
\end{itemize}

\noindent In Figure \ref{F-asint-cuartica}, we plot the curve
$\cal C$, and the asymptotes $\widetilde{\cal C}_1$ and $\widetilde{\cal C}_2$.
%\vspace*{-0.5cm}
\begin{figure}[h]
$$
\begin{array}{lcr}
\psfig{figure=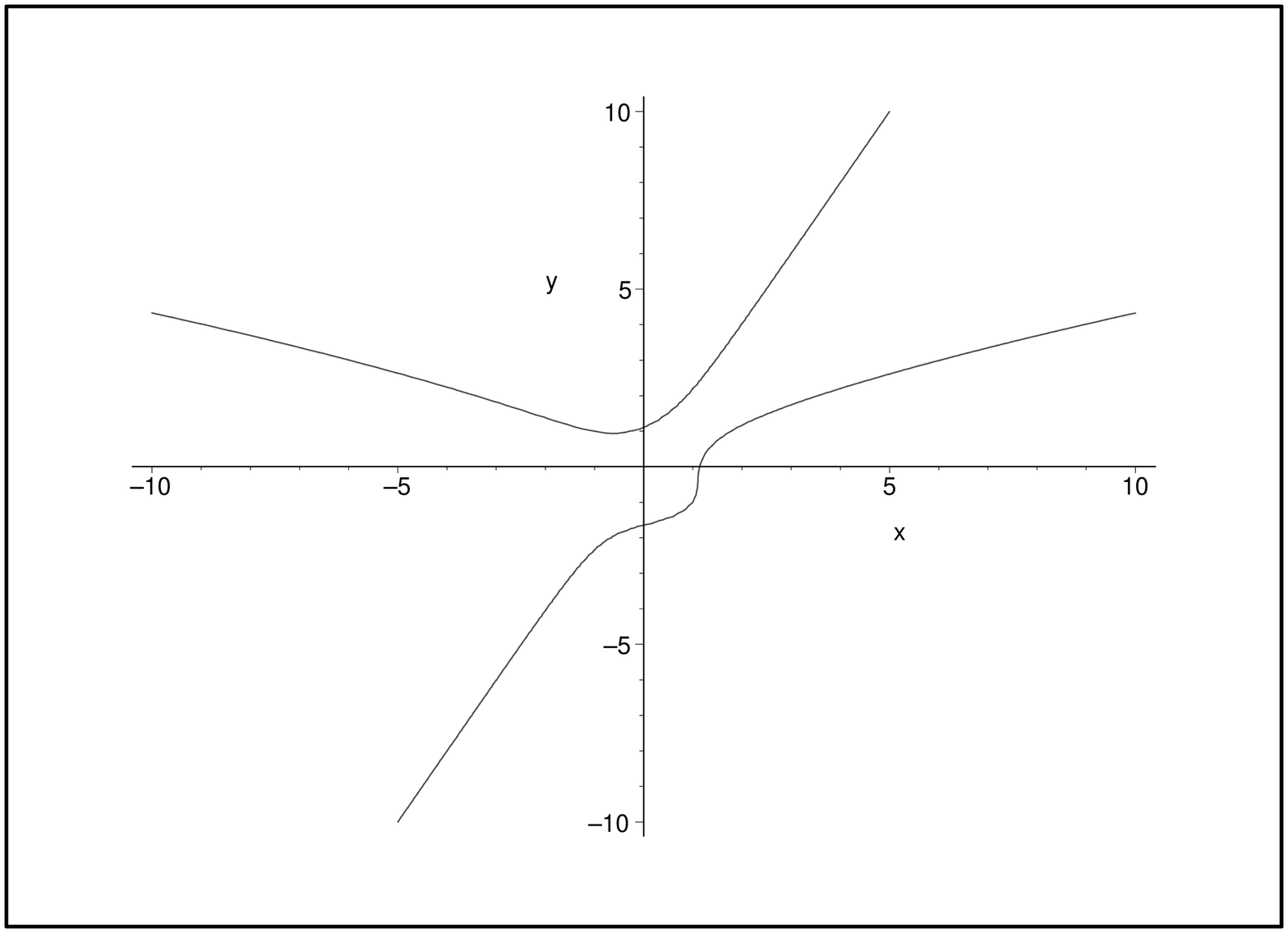,width=4.1cm,height=4.2cm,angle=270} &
\psfig{figure=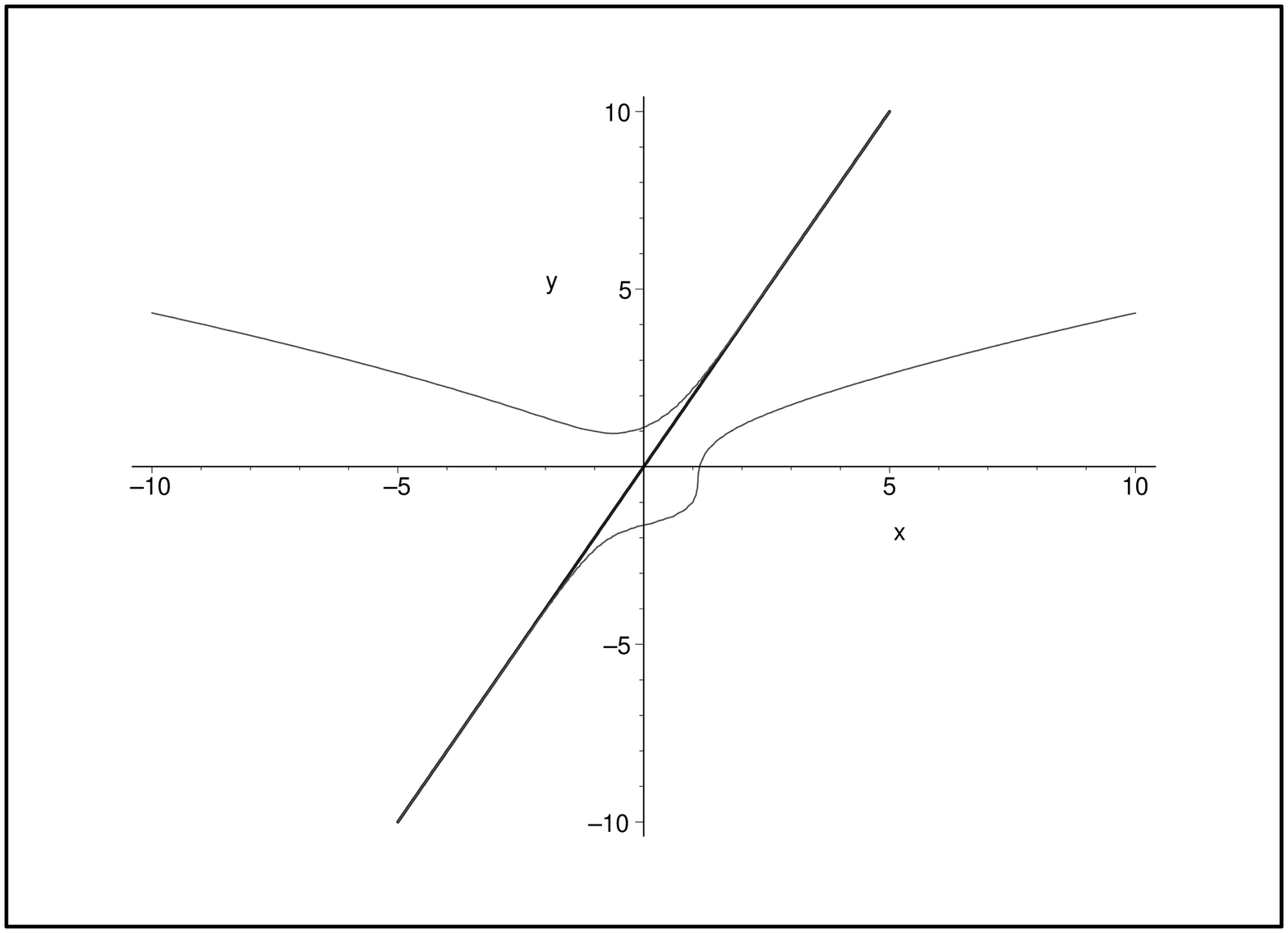,width=4.1cm,height=4.2cm,angle=270}&
\psfig{figure=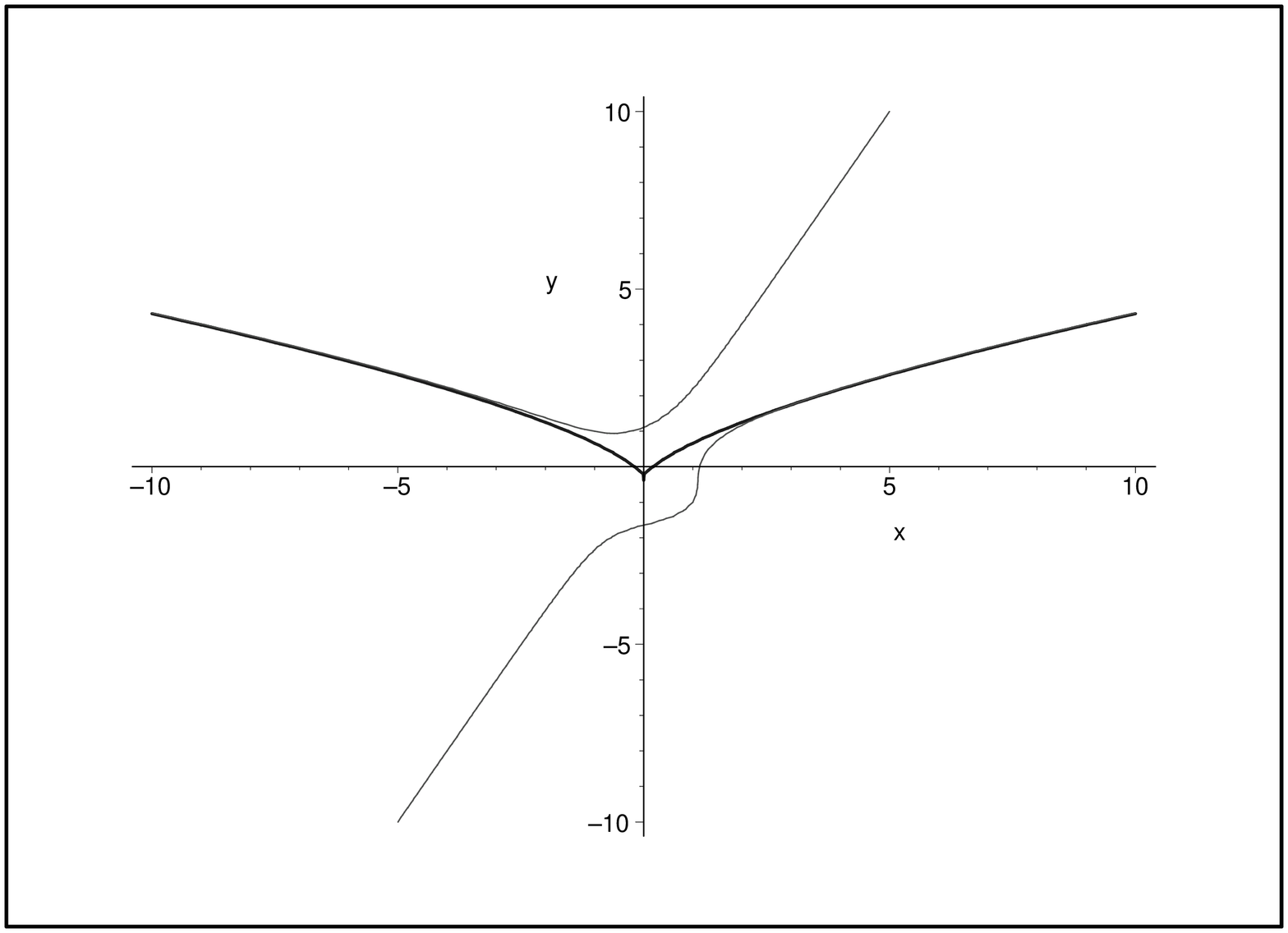,width=4.1cm,height=4.2cm,angle=270}
\end{array}
%\vspace*{-0.4cm}
$$ \caption{Curve $\cal C$ (left),  asymptote $\widetilde{\cal C}_1$ (center) and asymptote $\widetilde{\cal C}_2$ (right).}\label{F-asint-cuartica}
\end{figure}

\end{example}

\section{Some results on perfect curves}\label{S-perfect-curves}

An indispensably tool in the development of the results presented in this paper is
the notion  of perfect curve. In this section we try to understand better
this concept. For this purpose, in the following we show some properties concerning perfect curves.   We start with a necessary condition for a curve to
be perfect.

\begin{proposition}\label{P-perfect-necessary}  A perfect curve is polynomial (i.e. it admits a polynomial parametrization).
%Every perfect curve admits a polynomial parametrization.
\end{proposition}

\noindent\textbf{Proof:} Taking into account the results in \cite{Manocha}, one has that a curve is polynomial if and only if it has
 only one infinity branch. Thus, in order to prove
the proposition,  we only need to show that a perfect curve cannot
have more than one infinity branch.

\para

 For this purpose, let $\cal C$ be a perfect curve
 defined by the polynomial $f(x,y)$, and let us assume
that $\cal C$ has two different infinity branches,
$B_j=\{(z,r_j(z))\in {\Bbb C}^2:\,z\in {\Bbb C},\,|z|>M_j\},\,\,
j=1,2$, where  $r_j$ are as in equation (\ref{Eq-inf-branchN}), i.e.
$$r_1=m_1z+a_1z^{-\frac{r_1}{N_1}+1}+a_2z^{-\frac{r_2}{N_1}+1}+\cdots,\,\, r_2=m_2z+b_1z^{-\frac{s_1}{N_2}+1}+b_2z^{-\frac{s_2}{N_2}+1}+\cdots,$$
and $\nu(B_j)=N_j,\,j=1,2$  (see Remark \ref{R-conjugation}).

\para

\noindent Now, reasoning as in the proof of Lemma \ref{L-case3}, we
consider
$$h_1(x,y)=\prod_{j=1}^{N_1}(y-r_{1,j}(x))\quad\quad h_2(x,y)=\prod_{j=1}^{N_2}(y-r_{2,j}(x))$$
where $r_{i,1}(z),\ldots r_{i,N_i}(z)$ are the conjugates
of   $r_i(z)$,\,$i=1,2$.

\para

It holds that $h_1$ and $h_2$ divide $f$ w.r.t. the variable $y$ (see
proof of Lemma \ref{L-case3}). In addition, $h_1$ and $h_2$ do not have common
factors since  $r_{1,i}$ and $r_{2,j}$ belong to different conjugacy classes for every $i\in \{1,\ldots,N_1\}$, and $j\in \{1,\ldots,N_2\}$.
Thus, we deduce that
 $$\deg({\cal C})\geq \deg_y(f)\geq N_1+N_2> N_1\geq n_1,\quad \mbox{where $\deg(B_1)=n_1$.}$$
 In addition, from Theorem
\ref{T-constr-asintota}, we have that $\cal C$ can be approached at
$B_1$ by an asymptote of degree $n_1$. Therefore, we get that
$\cal C$ is not perfect which contradicts the assumption. Hence, we
conclude that a perfect curve cannot have more than one infinity
branch. \hfill $\Box$

\para

In the following example, we show that  the reciprocal of
Proposition \ref{P-perfect-necessary} is not true. More precisely, we
consider a polynomial curve $\cal C$,  and we show that $\cal C$ is
not a perfect curve.

\begin{example}
Let $\cal C$ be the curve defined by the polynomial parametrization
$${\cal P}(t)=(t^4+t,t^2)\in {\Bbb R}[t]^2.$$ We compute the  implicit polynomial of $\cal C$  by
applying for instance the results in \cite{SWP} (see Chapter 4). We have that
$$f(x,y)=\resultant_t(x-t^4-t,y-t^2)=-y+x^2-2xy^2+y^4\in {\Bbb R}[x,y].$$
We apply Algorithm
{\sf Asymptotes Construction} to determine the asymptotes of $\cal C$. We first observe that $\cal C$ only has the infinity point  $P=(1:0:0)$. We compute the infinity branch associated to $P$ (we use the {\sf algcurves} package included in {\sf Maple}), and we get that
$B=\{(z,r(z))\in {\Bbb C}^2:\,z\in {\Bbb C},\,|z|>M\},$ where
$$r(z)=z^{1/2}+\frac{1}{2}z^{-1/4}-\frac{1}{64}z^{-7/4}+\frac{1}{128}z^{-10/4}+\cdots.$$
Now, we consider
$\tilde{r}(z)=z^{1/2}$, and we obtain the asymptote $\widetilde{\cal C}$  defined by the parametrization
$\widetilde{Q}(t)=(t^2,t).$ The polynomial defining implicitly  $\widetilde{\cal C}$ is $\tilde{f}(x,y)=y^2-x$.

\para

\noindent Thus, the curve $\cal C$ that has degree 4, is
approached by  $\widetilde{\cal C}$ that has degree 2. Therefore, $\cal C$ is not perfect. In
Figure \ref{F-polinomial-noperfecta}, we plot  $\cal
C$ and the asymptote $\widetilde{\cal C}$.
%\vspace*{-0.6cm}
\begin{figure}[h]
$$
\begin{array}{lr}
\psfig{figure=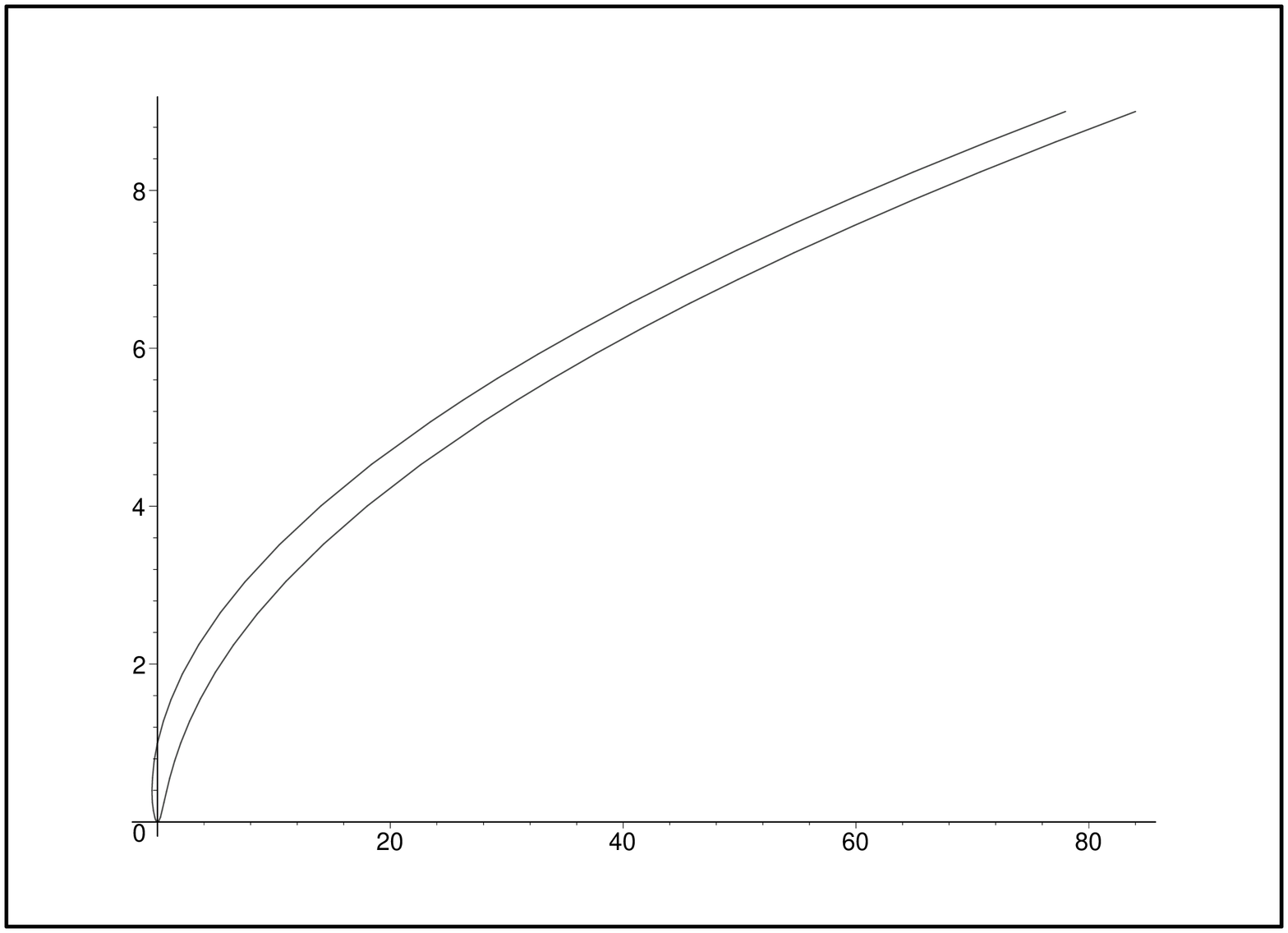,width=4cm,height=4cm,angle=270} &
\psfig{figure=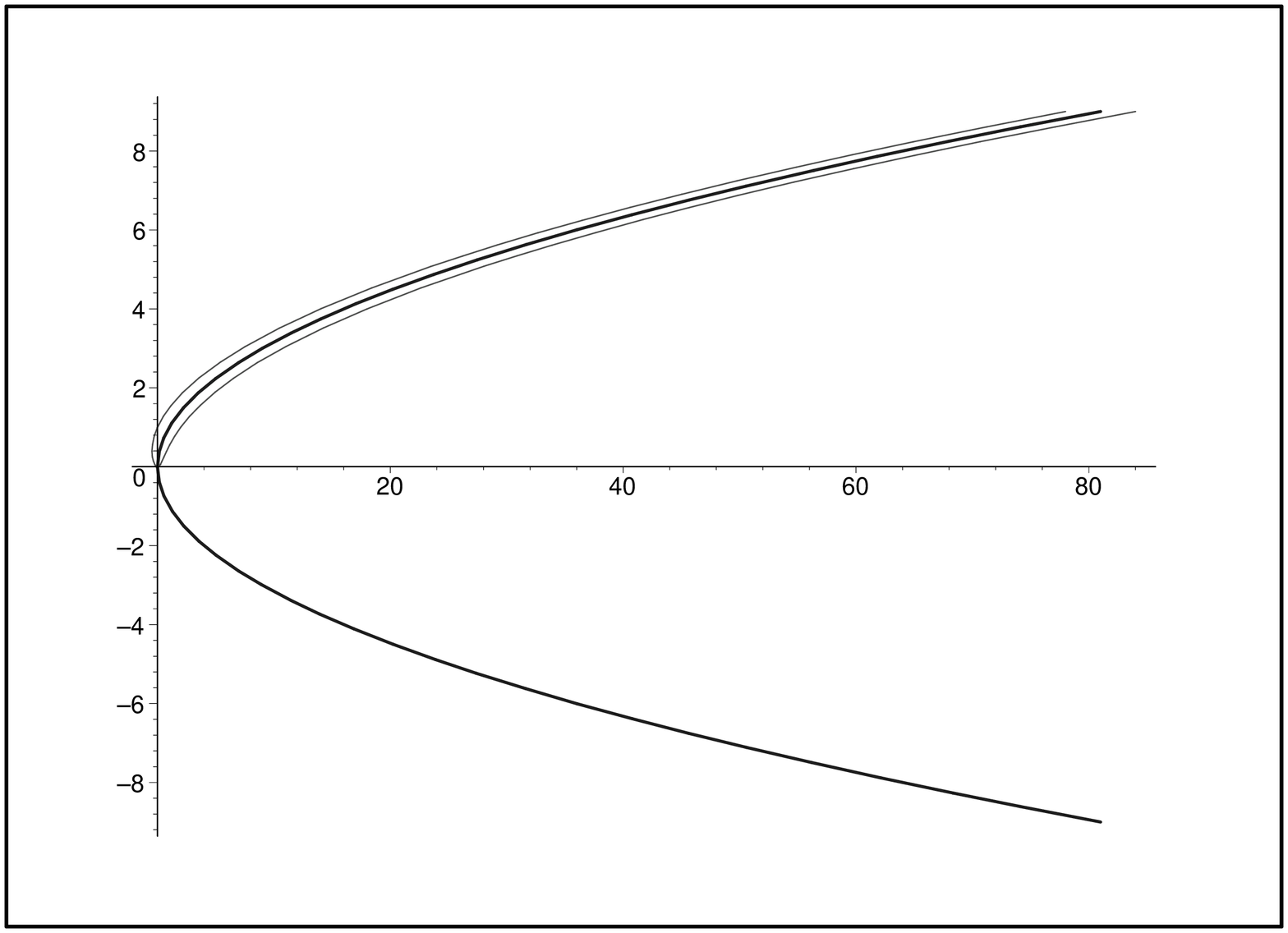,width=4cm,height=4cm,angle=270}
\end{array}
\vspace*{-0.4cm}
$$ \caption{Curve $\cal C$ (left), and asymptote $\widetilde{\cal C}$ (right).}\label{F-polinomial-noperfecta}
\end{figure}

\para

We remark that  $\widetilde{\cal C}$   approaches $\cal C$    at every leaf of the branch $B$ (see statement 3 in Remark \ref{R-approaching-curves}). More precisely, $B$ has two real and two complex leaves (see Remark 4.7 in
\cite{paper1}). The real leaves are convergent with the leaf
$\tilde{r}_1(z)=z^{1/2}$ of the parabola $\widetilde{\cal C}$, and the complex leaves are convergent with the leaf
$\tilde{r}_2(z)=-z^{1/2}$ of $\widetilde{\cal C}$.
% \textcolor{red}{At this point, note that $B$ has two complex
%leafs %(see Definition \ref{D-infinitybranch})
%$$r_3(z)=-z^{1/2}+\frac{I}{2}z^{-1/4}+\frac{I}{64}z^{-7/4}-\frac{1}{128}z^{-10/4}+\cdots$$
%$$r_4(z)=-z^{1/2}-\frac{I}{2}z^{-1/4}-\frac{I}{64}z^{-7/4}-\frac{1}{128}z^{-10/4}+\cdots$$
%(which cannot be plotted) convergent with the leaf
%$\tilde{r}_2(z)=-z^{1/2}$ of the parabola .}
\end{example}

\para

Proposition \ref{P-perfect-necessary} states that a perfect curve has  only one infinity
branch but this condition does not ensure that the curve is perfect. In the following, we provide  a characterization of perfect curves (see Proposition \ref{P-perfect-characterization}). For this purpose, we first remark some properties obtained from the definitions and results introduced before.

\begin{remark}\label{R-same-degree}
\begin{enumerate}
\item Two approaching perfect curves   have the same degree (see Definition \ref{D-perfect-curve}).
\item Two asymptotes that approach the same branch  have the same degree (see Corollary \ref{C-approaching-curves}).
\item Any asymptote that approaches a curve at a branch of degree $n$  has degree $n$ (see  Remark \ref{R-grado-caso1} and Theorem \ref{T-constr-asintota}).
\item From statement above and Lemma \ref{L-case3}, we deduce that if $B$ is a branch of a perfect curve then $\nu(B)=\deg(B)$.
\end{enumerate}
\end{remark}

\begin{proposition}\label{P-perfect-characterization}
Let $\cal C$ be an algebraic plane curve  of
degree $n$. $\cal C$ is perfect if and only if it has a unique
infinity branch $B$, and $\deg(B)=n$.
\end{proposition}

\noindent\textbf{Proof:} First, we assume that $\cal C$ is perfect.   From
Proposition \ref{P-perfect-necessary}, we get that $\cal C$ only has one
infinity branch,  $B$. It holds that $\deg(B)\geq n$; otherwise $\cal C$ would be approached by asymptotes of degree less
than $n$ (see Remark \ref{R-same-degree}, statement 3). Moreover, from Remark
\ref{R-grado-caso3}, we get that $\deg(B)\leq n$. Thus, we deduce that $\deg(B)=n$.\\
Now, let $\cal C$ be such that it has a unique infinity branch $B$,
and $\deg(B)=n$. Let us assume   that $\cal C$ is not perfect. Then,
there exists a curve $\overline{\cal C}$, with $\deg(\overline{\cal
C})<n$, that approaches $\cal C$ at  $B$. From
Remark \ref{R-minimal-degree}, we get that the
asymptotes of $\cal C$ at $B$ have degree less than $n$, which
contradicts Remark \ref{R-same-degree} (statement $3$). Therefore, we conclude
that $\cal C$ is perfect.\hfill $\Box$

\para

\begin{example}\label{Ej-YnXm}
Let $\cal C$ be the plane curve   defined by the irreducible polynomial $f(x,y)=y^n-x^m\in {\Bbb R}[x,y],\,n,m\in {\Bbb N},\,n>m$. Let us prove that $\cal C$ is perfect. For this purpose, we first observe that  $\gcd(n,m)=1$; otherwise,
$n=n_1k,\, m=m_1k,\,k\geq 2$, and $(y^{n_1}-x^{m_1})$ divides $f$
%$$y^n-x^m=(y^{n_1}-x^{m_1})(y^{n_1(k-1)}+y^{n_1(k-2)}x^{m_1}+\cdots+y^{n_1}x^{m_1(k-2)}+x^{m_1(k-1)})$$
which is impossible because $f$ is irreducible.
%In addition, we assume w.l.o.g that $n>m$, and then $\deg({\cal C})=n$.
%otherwise, the unique infinity point of $\cal C$ is $P=(0:1:0)$,
%\textcolor{red}{and we have discarded this case at the beginning of
%the paper}. As a consequence, $\cal C$ has degree $n$.

\para

Note that $\cal C$ has a unique infinity branch $B$, since it admits the
polynomial parametrization $(t^n,t^m)$. In addition, $B$ is given by $r(z)=z^{m/n}$, and thus $\deg(B)=n$.
Then, by applying Proposition \ref{P-perfect-characterization}, we
conclude that $\cal C$ is perfect.
\end{example}

\para

The following proposition states a sufficient condition for a curve
to be perfect, which can be checked without computing any infinity
branch. This result will play an important role in Section
\ref{S-families}.

\begin{proposition}\label{P-perfect-sufficient}
Let $\cal C$ be an algebraic plane curve with a unique infinity
point $P$, and let $P$ be regular. Then $\cal C$ is perfect.
\end{proposition}

\noindent\textbf{Proof:} Let $\degree({\cal C})=d$, and we assume
that $d\geq 2$ (the case of lines is trivial). Let us prove that
${\cal C}$ satisfies the conditions of Proposition
\ref{P-perfect-characterization}; i.e, ${\cal C}$ has a unique
infinity branch, say $B$, and $\deg(B)=d$.

\para

We assume w.l.o.g that
$P=(1:0:0)$ (otherwise, we apply a linear change of coordinates).
Since $P$ is regular, there exists a unique place centered at $P$
and hence, there exists a unique infinity branch associated to $P$
(see Section 2). Therefore the first condition
holds.

\para

Now, we focus on the second condition. Let $B=\{(z,r(z))\in {\Bbb
C}^2: \,z\in {\Bbb C},\,|z|>M\}$,  $M\in {\Bbb R}^+$, be the
unique infinity branch of $\cal C$, where
\[r(z)=a_1z^{-N_1/N+1}+\cdots
+a_kz^{-N_k/N+1}+a_{k+1}z^{-N_{k+1}/N+1}+\cdots,\]
$a_1,a_2,\ldots\in\mathbb{C}\setminus \{0\}$,
$N, N_1,N_2\ldots\in\mathbb{N}$,  $0<N_1<N_2<\cdots$, and
  $N_k\leq N <N_{k+1}$  (see equation (\ref{Eq-inf-branchN})). Let $\gcd(N,N_1,\ldots,N_k)=b$,\,$N_j=n_jb,\,N=nb,\,\,j=1,\ldots,k$ (see equation (\ref{Eq-inf-branchn})); that is,
  $\deg(B)=n$ (see Definition \ref{D-degreebranch}). In the following, we prove that $n=d$. For this purpose, we analyze how $r(z)$ is obtained. Let
$$f(x,y)=y^d+f_{d-1}(x,y)+f_{d-2}(x,y)+\cdots+f_1(x,y)+f_0,$$
where $f_i(x,y)$ is the homogeneous form of degree $i$. In addition,
let
$$f_{d-1}(x,y)=b_0x^{d-1}+b_1x^{d-2}y+\cdots+b_{d-2}xy^{d-2}+b_{d-1}y^{d-1},\quad \mbox{and}$$
$$f_{d-2}(x,y)=c_0x^{d-2}+c_1x^{d-3}y+\cdots+c_{d-3}xy^{d-3}+c_{d-2}y^{d-2}.$$
The projective curve associated to $\cal C$ is given by:
$$F(x:y:z)=y^d+zf_{d-1}(x,y)+z^2f_{d-2}(x,y)+\cdots+z^{d-1}f_1(x,y)+z^df_0.$$
We observe that  $b_0\not=0$. Indeed: since $d\geq 2$, we note that
$$\frac{\partial F}{\partial y}(P)=\frac{\partial f_d}{\partial
y}(1,0)=0,\quad\text{ and } \quad\frac{\partial F}{\partial
z}(P)=f_{d-1}(1,0)=b_0.$$ Furthermore, from Euler's formula, we have
that $\frac{\partial F}{\partial x}(P)=0$. Then, since  $P$ is
regular, we deduce that $b_0\neq 0$.

\para

Under these conditions, in order to compute $r(z)$, we consider the polynomial (see Section 2):
$$g(y,z)=F(1:y:z)=y^d+zf_{d-1}(1,y)+\cdots+z^df_0$$
which can be written as
$$g(y,z)=y^d+z(b_0+b_1y+\cdots+b_{d-2}y^{d-2}+b_{d-1}y^{d-1})+$$
$$z^2(c_0+c_1y+\cdots+c_{d-3}y^{d-3}+c_{d-2}y^{d-2})+\cdots.$$
In addition, we have that $r(z)=z\varphi(z^{-1})$, where
$\varphi(z)$ is a series expansion for a solution  of
$g(y,z)=F(1:y:z)=0$; that is, $g(\varphi(z), z)=0$, where
$\varphi(z)=a_1z^{N_1/N}+a_2z^{N_2/N}+a_3z^{N_3/N} + \cdots$ (see
Section \ref{S-notation}). Hence, the expression
$$g(\varphi(z),z)=(a_1z^{N_1/N}+a_2z^{N_2/N}+\cdots)^d+$$
$$z\left(b_0+b_1(a_1z^{N_1/N}+a_2z^{N_2/N}+\cdots)+b_2(a_1z^{N_1/N}+a_2z^{N_2/N}+\cdots)^2+\cdots\right)+$$
$$z^2\left(c_0+c_1(a_1z^{N_1/N}+a_2z^{N_2/N}+\cdots)+c_2(a_1z^{N_1/N}+a_2z^{N_2/N}+\cdots)^2+\cdots\right)+\cdots$$
vanishes for $|z|<M^{-1}$ which implies that terms
with a common exponent must cancel.

\para

Since  $b_0\neq 0$, the terms with lowest order in $g(\varphi(z),z)$
are $b_0z$ and $(a_1z^{N_1/N})^d=a_1^dz^{N_1d/N}$. Thus, they must
cancel and then,  $N=N_1d$ and $n=n_1d$ (note that $N=nb$, and $N_1=n_1b$). Hence,  since $n_1\geq 1$, we get that
$n\geq d$. On the other hand, Remark \ref{R-grado-caso3} states that
$d\geq N \geq n$, so $n=d$.

\para

Therefore,   $\cal C$ has a unique infinity
branch $B$, and   $\deg(B)=d$. From Proposition
\ref{P-perfect-characterization}, we conclude that $\cal C$ is a perfect
curve. \hfill $\Box$

\para

\begin{remark}
The reciprocal of Proposition \ref{P-perfect-sufficient} is not
true. For instance, the curve $y^3-x=0$ is perfect (see Example
\ref{Ej-YnXm}) but its unique infinity point, $(1:0:0)$, is singular.
\end{remark}

\section{Families of asymptotes}\label{S-families}

In Section \ref{S-asymptotes}, given an algebraic plane
curve $\cal C$, and  an infinity branch $B$ of $\cal C$, we prove that $\cal C$ can be approached by an
asymptote at $B$. In the following, we  see that this asymptote may not be unique. In fact,   in most of cases, there are infinitely many asymptotes
associated to a given infinity branch.

\para

\begin{figure}[h]
\begin{center}
\psfig{figure=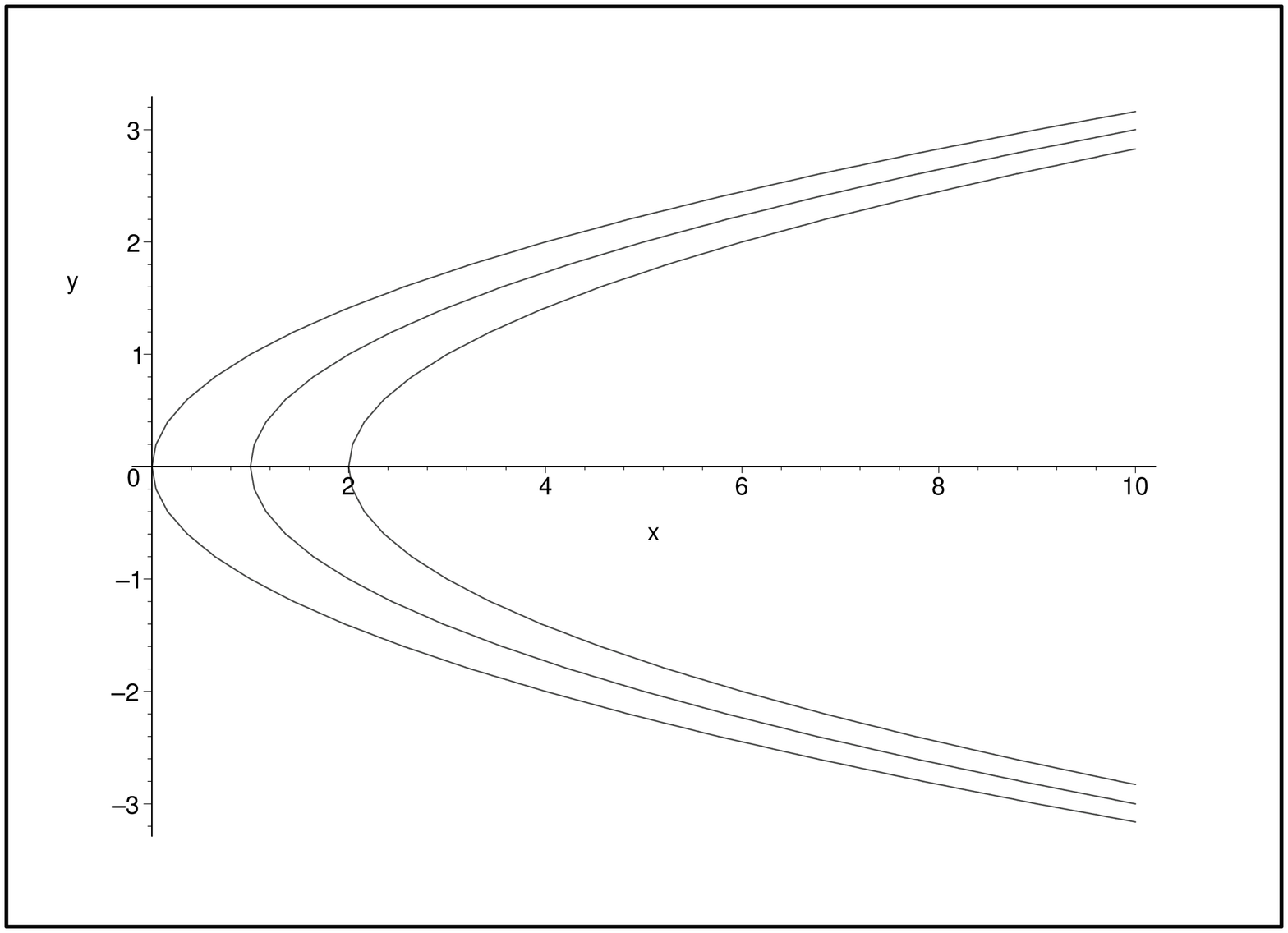,width=4cm,height=6cm,angle=270}
\end{center}
\vspace*{-0.5cm}
\caption{Curves ${\cal C}$, ${\cal D}_1$ and ${\cal
D}_2$.}\label{F-proximity-class}
\end{figure}

For instance, let $\cal C$ be the curve defined by the polynomial
$f(x,y)=y^2-x$ and, for each $k\in\mathbb{R}$, let ${\cal D}_k$ be
the curve defined by $g_k(x,y)=y^2-x+k$. It holds that  ${\cal D}_k$
approaches $\cal C$, for every $k\in\mathbb{R}$. Indeed: note that
both curves admit a polynomial parametrization, so each of them has
only one infinity branch. Furthermore, these branches are defined by
$r(z)={z}^{1/2}$, and
$r_k(z)=z^{1/2}-k/(2z^{1/2})-k^2/(8z^{3/2})+\cdots=(z-k)^{1/2}$, %-k^3/(16z^{5/2})
respectively.
%(-k)^{1/2}(-1+z/(2k)+z^2/(8k^2)+z^3/(16k^3)+\cdots)
Clearly, these branches are convergent, since
$\lim_{z\rightarrow\infty}(r(z)-r_k(z))=0.$ Thus, $\{{\cal
D}_k\}_{k\in\mathbb{R}}$ defines  an infinite family of  curves that
approach  $\cal C$ at its unique infinity branch (see Figure
\ref{F-proximity-class}).

\para

Note that the relation of {\it proximity}   is actually an equivalence relation for the set of curves having only one infinity
branch. Clearly, it is reflexive, since every curve approaches
itself. Moreover, statement 1 in Remark \ref{R-approaching-curves} ensures the
symmetry.

\para

Transitivity does not hold for a general set of curves. More precisely, let
  ${\cal C}_1$, ${\cal C}_2$ and ${\cal C}_3$ be three plane curves such that
${\cal C}_1$ approaches ${\cal C}_2$ at an infinity branch $B$, and
  ${\cal C}_2$ approaches ${\cal C}_3$ at a different infinity
branch $B^*$. We can not get that ${\cal C}_1$ approaches ${\cal
C}_3$. However, if we consider curves with  only one
infinity branch, $B$ and $B^*$ are the same and  then, ${\cal
C}_1$ approaches ${\cal C}_3$ at $B$. In this case, transitivity holds and then,
we have an equivalence relation for the set of curves having only one infinity
branch. In the following, we refer to the equivalence classes
associated to this equivalence relation as
{\it proximity classes}.

\para

We observe that this property can be applied to perfect curves
since, from Proposition \ref{P-perfect-necessary}, perfect curves  only have one infinity branch. That is, we consider the above property  restricted to  the set of perfect curves and then,  given a perfect curve $\cal C$,
the set of perfect curves approaching $\cal C$ determines a proximity
class. Therefore, the set of asymptotes that approach a curve
at a given infinity branch is also a proximity class.

\para

 In this
section, we provide a method  that allows, in same
cases, to compute the curves in a proximity class. For this
purpose, we first introduce the following definition.

\begin{definition}\label{D-regular-perfect-curve}
A \emph{regular perfect curve} is a  curve having a unique infinity
point, which is regular.
\end{definition}

From Proposition \ref{P-perfect-sufficient}, regular perfect curves are a subset of the set of perfect curves. In the
following, we prove a nice property of this subset which will allow us
to decide whether two regular perfect curves approach each other  by
comparing their implicit polynomials.

\para

More precisely, Theorem \ref{T-irrel-terms} states that a
polynomial defining implicitly a regular
perfect curve has some {\it irrelevant terms} (i.e. terms that
do not affect   the asymptotic behavior of the curve). If we
modify the coefficients of these terms, we obtain a new perfect curve
that approaches the original one. Hence, we can construct infinitely many
different asymptotes approaching a given plane curve simply by manipulating its
irrelevant terms.

\begin{theorem}\label{T-irrel-terms}
Let $\cal C$ be a regular perfect curve  defined by an irreducible
polynomial $f\in\mathbb{R}[x,y]$ of degree $d$. The asymptotic
behavior of $\cal C$ is completely determined by $f_d$ and $f_{d-1}$. %the terms of degree $d$ and $d-1$ in $f$
\end{theorem}

\noindent
\textbf{Proof:} Taking into account Proposition
\ref{P-perfect-necessary}, we have that $\cal C$  has  only one
infinity branch. Let  $B=\{(z,r(z))\in {\Bbb C}^2:
\,z\in {\Bbb C},\,|z|>M\}$ be this branch. From Proposition \ref{P-perfect-characterization}, we get that
$\deg(B)=d$ which implies that $r(z)$ can be written as
$$r(z)=mz+a_1z^{1-1/d}+a_2z^{1-2/d}+\cdots +a_{d-1}z^{1-(d-1)/d}+a_d+A(z)$$
where $m,a_1,a_2,\ldots\in\mathbb{C}$, and $A(z)=\sum_{j=d+1}^\infty
a_{j}z^{-q_j}$, $q_j\in\mathbb{Q}^+$ for $j\geq d+1$.
Note that this notation is slightly different from
that of (\ref{Eq-inf-branchn}); here, some coefficients $a_i$ might
be zero.

\para

In order to describe the asymptotic behavior of $\cal C$, we just
need to compute the terms with non negative exponent in
$r(z)$. That is, we need to determine the coefficients
$m,a_1,\ldots a_d$. The first of them, $m$, can be obtained
directly from the homogeneous form of maximum degree,
$f_d(x,y)=(y-mx)^d$. In the following we assume w.l.o.g. that $m=0$; that is, $B$ is associated to the infinity point $P=(1:0:0)$ (otherwise, we apply a linear change of coordinates).

\para

The  coefficients $a_1,a_2,\ldots a_d$ can be obtained by applying the
procedure described in the proof of Proposition
\ref{P-perfect-sufficient}. The condition to be held is that
$g(\varphi(z),z)=0$ for $|z|<M^{-1}$, where
$$\varphi(z)=zr(z^{-1})=a_1z^{1/d}+a_2z^{2/d}+\cdots +a_{d-1}z^{(d-1)/d}+a_dz+zA(z^{-1}).$$
Reasoning as in the proof of Proposition
\ref{P-perfect-sufficient}, we get that  the terms with a common exponent in the expression
$$g(\varphi(z),z)=(a_1z^{1/d}+a_2z^{2/d}+\cdots)^d+$$
$$+z\left(b_0+b_1(a_1z^{1/d}+a_2z^{2/d}+\cdots)+b_2(a_1z^{1/d}+a_2z^{2/d}+\cdots)^2+\cdots\right)+$$
$$+z^2\left(c_0+c_1(a_1z^{1/d}+a_2z^{2/d}+\cdots)+c_2(a_1z^{1/d}+a_2z^{2/d}+\cdots)^2+\cdots\right)+\cdots$$
must cancel. In addition, we also have that $b_0\not=0$.

\para

First, we consider the terms with minimum exponent. They
are $b_0z$ and $a_1^dz$. Then, $b_0+a_1^d=0$
which implies that $a_1=(-b_0)^{1/d}$. We set $\textsf{A}_1:=a_1$ ($\textsf{A}_1$  represents
that $a_1$ has now a fixed value).

\para

We substitute the value of $a_1$ in $g(\varphi(z),z)$, and we reason similarly to compute
$a_2$. Note that the terms $b_0z$ and
$a_1^dz$ have been canceled, and $a_1$ has a fixed value $\textsf{A}_1\not=0$ (note that  $b_0\not=0$).
Hence, the terms with minimum exponent are
$b_1\textsf{A}_1z^{1+\frac{1}{d}}$ and $d\textsf{A}_1^{d-1}a_2z^{1+\frac{1}{d}}$. Then, $b_1\textsf{A}_1+d\textsf{A}_1^{d-1}a_2=0$ which implies that the value of $a_2$ is determined  from $b_{1}$.  We have  that $a_2=\frac{-b_1}{d\textsf{A}_1^{d-2}}$.  We set $\textsf{A}_2:=a_2$.

\para

%For better understanding the process, we consider one more step.
Once we have $a_1$ and $a_2$,  we compute $a_3$. For this purpose, we substitute  the values of $a_1, a_2$ in $g(\varphi(z),z)$, and the terms with
minimum exponent are $b_2\textsf{A}_1^2$, $b_1\textsf{A}_2$,
$\left(\begin{array}{c}d\\2\end{array}\right)\textsf{A}_1^{d-2}\textsf{A}_2^2$ and
$d\textsf{A}_1^{d-1}a_3$ (all of them are multiplied by $z^{1+\frac{2}{d}}$). Thus, we
have that
$$b_2\textsf{A}_1^2+b_1\textsf{A}_2+\left(\begin{array}{c}d\\2\end{array}\right)\textsf{A}_1^{d-2}\textsf{A}_2^2+d\textsf{A}_1^{d-1}a_3=0.$$
Again, we obtain an equation where $a_3$, $b_2$ and a set of constants appear. Then, the value of $a_3$ is determined from $b_{2}$ (we recall that $\textsf{A}_1\not=0$).

\para

 In the following, we  prove that once the values of $a_j,\,j=1,\ldots, i-1$ are obtained, reasoning as above, we   get  an equation where   $a_i$, $b_{i-1}$ and a set of constants appear. %Then,  the value of $a_i$ is determined from $b_{i-1}$. %this is a general rule.
For this purpose, we first observe that the term with minimum exponent including  $a_i$  is $d\textsf{A}_1^{d-1}a_iz^{1+\frac{i-1}{d}}$, and
the term with minimum exponent including $b_{i-1}$ is
$b_{i-1}\textsf{A}_1^{i-1}z^{1+\frac{i-1}{d}}$ (note that  $\textsf{A}_1\not=0$). Hence, both coefficients,
$a_i$ and $b_{i-1}$, appear in the equation multiplied by $z^{1+\frac{i-1}{d}}$.
%Therefore, both coefficients,
%$a_i$ and $b_{i-1}$, enter the equation, for the first time, in the
%$i$-th step, both multiplied by $z^{1+\frac{i-1}{d}}$.
The remainder elements involved in this equation are constants computed in
the previous steps. Thus,    $a_i$
 is obtained from $b_{i-1}$.

\para

Therefore, we conclude that the asymptotic behavior of the curve is
determined by the coefficient $m$, obtained from $f_{d}(x,y)$, and
the coefficients $b_0,b_1,\ldots,b_{d-1}$, obtained from
$f_{d-1}(x,y)$.\hfill $\Box$

\para

Theorem \ref{T-irrel-terms} implies that any modification in the terms
of $f$ of degree less or equal to $d-2$, does not affect the asymptotic
behavior of the curve (we refer to them as the {\it irrelevant terms}). In addition, from the proof of Theorem \ref{T-irrel-terms}, we also get that
%However, we have just proved something more. We have shown that, for
for regular perfect curves, there exists a one to one correspondence
between the coefficients $m,a_1,a_2,\ldots a_d$ (that determine the
asymptotic behavior of the curve), and the coefficients of
$f_{d}$ and $f_{d-1}$. Hence, we deduce the following
corollary.

\begin{corollary}\label{C-irrel-terms}
Two regular perfect curves approach each other if and only if their
terms of degree $d$ and $d-1$ are the same.
\end{corollary}

 In the following, we consider two approaching perfect curves ${\cal C}$ and $\overline{\cal C}$. From statement 2 in Remark \ref{R-approaching-curves}, we have that both curves   have
a unique common infinity point $P$. In the next corollary, we prove that   $P$ is a regular point of ${\cal C}$ iff $P$ is a regular point of $\overline{\cal C}$. Thus, we deduce that  a perfect curve having a singular infinity point
 cannot be approached by a regular perfect curve.

 \begin{corollary}\label{C-same-character} Let ${\cal C}$ and $\overline{\cal C}$ be two  approaching perfect curves. ${\cal C}$ is a regular perfect curve if and only if $\overline{\cal C}$ is a regular perfect curve.
\end{corollary}
\noindent\textbf{Proof:} Let $\cal C$ be defined  by the polynomial
$$f(x,y)=(y-mx)^d+f_{d-1}(x,y)+f_{d-2}(x,y)+\cdots+f_1(x,y)+f_0$$
where $f_i(x,y)$ is the homogeneous form of degree $i$, and $d\geq
2$ (the case of lines trivially holds). The projective curve
associated to $\cal C$ is given by
$$F(x:y:z)=(y-mx)^d+zf_{d-1}(x,y)+z^2f_{d-2}(x,y)+\cdots+z^{d-1}f_1(x,y)+z^df_0.$$
Thus, $\frac{\partial F}{\partial y}(P)=0$ ($d\geq 2$),
$\frac{\partial F}{\partial z}(P)=f_{d-1}(1,m),$ and $\frac{\partial
F}{\partial x}(P)=0$ (apply Euler's formula). Then,  $P$ is regular
if and only if $f_{d-1}(1,m)\not=0$. Finally, the result follows by
applying Corollary \ref{C-irrel-terms}.\hfill $\Box$

\para
\para

Given a regular perfect curve $\cal C$ of degree $d$
defined by a polynomial $f(x,y)$,  we
may compute all the curves in
its proximity class simply by  modifying the irrelevant terms in $f$ (see Corollaries \ref{C-irrel-terms} and \ref{C-same-character}).

\para

\noindent
For instance, let $\cal C$ be the curve defined by the polynomial
$$f(x,y)=x^3+3x^2y+3xy^2+y^3+2x^2+y-3\in {\Bbb R}[x,y].$$
$\cal C$ has   only one
infinity point, $P=(1:-1:0)$, which is regular. Thus,  $\cal C$ is a regular
perfect curve. The curves within the
  proximity class   of $\cal C$ are
 implicitly defined by the polynomials
$$x^3+3x^2y+3xy^2+y^3+2x^2+ax+by+c,\quad a,b,c\in\mathbb{R}. $$
 Note that any curve that belongs
to this proximity class can be associated  uniquely to a
vector $(a,b,c)\in\mathbb{R}^3$. This remark motivates the following
proposition.

\begin{proposition} \label{P-evectorial}  Let $\cal C$ be a regular perfect curve of degree $d$. The proximity class of $\cal C$ is isomorphic to $\mathbb{R}^{\frac{(d-1)d}{2}}$. \end{proposition}
%A proximity
%class of the  set of regular perfect curves is isomorphic to $\mathbb{R}^{\frac{(d-1)d}{2}}$. Indeed: from Remark \ref{R-same-degree}, statement 1, we have that all the
%curves in the same proximity class  have the same degree $d$. Thus, the result follows taking into account that   the number of irrelevant terms in
%a generic polynomial of degree $d$ is $\frac{(d-1)d}{2}$.
%------------------------
%For a class of
%degree $d$ (from Remark \ref{R-same-degree}, statement 1, all the
%curves in the same proximity class  have the same degree $d$), the
%dimension of $\cal V$ is given by the number of irrelevant terms in
%a generic polynomial of degree $d$. That is,
%$\dim({\cal V})=\frac{(d-1)d}{2}.$
\noindent\textbf{Proof:} From Remark \ref{R-same-degree}, statement 1, we have that all the
curves in the same proximity class  have the same degree $d$. Thus, the result follows taking into account that   the number of irrelevant terms in
a generic polynomial of degree $d$ is $\frac{(d-1)d}{2}$, and that any curve that belongs
to this proximity class can be associated  uniquely to a
vector in $\mathbb{R}^{\frac{(d-1)d}{2}}$. \hfill $\Box$

\para

%\textcolor{red}{Finally, we observe that Theorem \ref{T-irrel-terms} provides a method for computing all the
%asymptotes that approach a given curve $\cal C$ at an infinity branch.}\footnote{esto no es cierto si la asintota que obtenemos no es regular}
%\textcolor{red}{ More precisely, let
%$\cal C$ be an algebraic plane curve with an infinity branch $B$, and
%let $\widetilde{{\cal C}}$ be the asymptote of $\cal C$ at $B$
%obtained by applying algorithm {\sf Asymptotes Construction}. If
%$\widetilde{{\cal C}}$ is a  regular perfect curve, we  may compute all the
%curves in its proximity class by modifying its irrelevant terms. The
%elements of this class are asymptotes of $\cal C$ at $B$.}\footnote{son todas?}

In the following, we show that Theorem \ref{T-irrel-terms} provides a method for computing all the
asymptotes of a curve $\cal C$ at an infinity branch $B$ associated to a regular infinity point.

\begin{theorem} \label{T-regular-asintota} Let $\cal C$ be a  plane curve, and $P$ a regular infinity point of $\cal C$. Let $B$ be an infinity branch
of $\cal C$ associated to  $P$, and  $\widetilde{{\cal C}}$  the asymptote of $\cal C$ at $B$
 obtained from algorithm {\sf Asymptotes Construction}. It holds that:
 \begin{itemize}\item[1.] $P$ is a regular point of $\widetilde{{\cal C}}$.
 \item[2.] The asymptotes of $\cal C$ at $B$ are the curves within the
  proximity class of  $\widetilde{{\cal C}}$.
 \end{itemize}
\end{theorem}
\noindent\textbf{Proof:}  Statement 2  is deduced
from Corollary \ref{C-irrel-terms}. So, in the following, we prove
statement 1. For this purpose, let $P=(1:0:0)$ (otherwise, we
consider a linear change of variables), and   $\deg(B)=n$.   Let
  $B=\{(z,r(z))\in {\Bbb C}^2:\,z\in {\Bbb C},\,|z|>M\}$,  where
$$r(z)=a_1z^{1-1/n}+a_2z^{1-2/n}+\cdots +a_{n-1}z^{1-(n-1)/n}+a_n+A(z)$$
$a_1,a_2,\ldots\in\mathbb{C}$, and $A(z)=\sum_{j=n+1}^\infty
a_{j}z^{-q_j}$, $q_j\in\mathbb{Q}^+$ for $j\geq n+1$ (note that this notation is slightly different from
that of (\ref{Eq-inf-branchn}); here, some coefficients $a_i$ might
be null).  From Algorithm {\sf Asymptotes Construction}, we get that a proper polynomial
  parametrization of $\widetilde{{\cal C}}$ is  given by $$\widetilde{Q}(t)=(q_1(t), q_2(t))=(t^{n},\,a_{1}t^{-1+n}+a_{2}t^{-2+n}+\cdots
+a_{n})\in {\Bbb C}[t]^2. $$
Under these conditions, we apply the results in \cite{Perez-Singularidades} and  Lemma \ref{L-case1}, and we have that the multiplicity of the point $P$ in $\widetilde{{\cal C}}$ is
$$\deg(\widetilde{{\cal C}})-{\displaystyle
\deg\left(\frac{1}{q_{2}(t)}\right)}=n-\deg(q_{2}).$$

Taking into account these previous results, we have to prove that
$\deg(q_{2})=n-1$, which is equivalent to $a_1\not=0$.

\para

For this purpose, since $P$ is an infinity point of $\cal C$, the
projective curve associated to $\cal C$ is given by
$$F(x:y:z)=f_d(x,y)+zf_{d-1}(x,y)+z^2f_{d-2}(x,y)+\cdots+z^{d-1}f_1(x,y)+z^df_0,$$
where
$$f_{d}(x,y)=y^s\prod_{j=1}^{d-s}(y-m_jx)=y^d+\ell_{d-1}y^{d-1}x+\cdots+\ell_{s}y^sx^{d-s},\,\,s\geq 1,$$
$$f_{d-1}(x,y)=b_0x^{d-1}+b_1x^{d-2}y+\cdots+b_{d-2}xy^{d-2}+b_{d-1}y^{d-1}.$$
%$$f_{d-2}(x,y)=c_0x^{d-2}+c_1x^{d-3}y+\cdots+c_{d-3}xy^{d-3}+c_{d-2}y^{d-2}.$$
If $s=1$, since $s\geq N\geq n$ (see Lemma
\ref{L-case3}), we get that $\deg(B)=n=1$ which implies that  $\widetilde{{\cal C}}$  is a
line. Thus, $P$ is a regular point in $\widetilde{{\cal C}}$. So,
let $s\geq 2$.  Then, $\frac{\partial F}{\partial x}(P)=0$ (apply
Euler's formula), $\frac{\partial F}{\partial
z}(P)=f_{d-1}(1,0)=b_0,$ and $\frac{\partial F}{\partial
y}(P)=\frac{\partial f_d}{\partial y}(1,0)=0$ ($s\geq 2$) which
implies that  $P$ is a regular point in ${{\cal C}}$ if and only if
$b_0\not=0$.\para

In the following, we prove that $b_0\not=0$ implies that $a_1\not=0$
and then, the result holds. For this purpose, we observe that  the
coefficients $a_1,a_2,\ldots a_n$ can be obtained by applying the
procedure described in the proof of Proposition
\ref{P-perfect-sufficient}. Since
$g(\varphi(z),z)=F(1,\varphi(z),z)=0$, where
$$\varphi(z)=zr(z^{-1})=a_1z^{1/n}+a_2z^{2/n}+\cdots +a_{n-1}z^{(n-1)/n}+a_nz+zA(z^{-1}),$$
reasoning as in the proof of Proposition
\ref{P-perfect-sufficient}, we get that  the terms with a common exponent in the expression $g(\varphi(z),z)=$
$$(a_1z^{1/n}+a_2z^{2/n}+\cdots)^d+\ell_{d-1}(a_1z^{1/n}+a_2z^{2/n}+\cdots)^{d-1}+\cdots+\ell_s (a_1z^{1/n}+a_2z^{2/n}+\cdots)^s$$
$$+z\left(b_0+b_1(a_1z^{1/n}+a_2z^{2/n}+\cdots)+b_2(a_1z^{1/n}+a_2z^{2/n}+\cdots)^2+\cdots\right)+\cdots$$
%$$+z^2\left(c_0+c_1(a_1z^{1/n}+a_2z^{2/n}+\cdots)+c_2(a_1z^{1/n}+a_2z^{2/n}+\cdots)^2+\cdots\right)+\cdots$$
must cancel.

\para

If $b_0\not=0$,  the first non-zero term with minimum exponent
appearing  is $b_0z$. Thus, there exists $j\in {\Bbb N}$ such that
$js/n=1$, and $b_0+a_j^s=0$, and $a_i=0,\,\,i=1,\ldots,j-1$. Since
$n\leq N\leq s$ (see Lemma \ref{L-case3}), we deduce that $j=1$,
 $n=s$, and $a_1=(-b_0)^{1/s}\not=0$, as we wanted
to prove.
%si $s=1$, entonces $n=1$ y la unica posibilidad que tenemos para que los terminos se vayan es que $a_1b_0\not=0$
\hfill $\Box$

\para

\begin{remark}
If $\widetilde{{\cal C}}$ is a line, there are no irrelevant terms. In this case the asymptote is unique.
\end{remark}

\noindent
In the following, we illustrate  these results with an example.

\begin{example}
We consider the curve $\cal C$ introduced in Example \ref{Ej-asint-cuartica}, and defined by the   polynomial
$$f(x,y)=2y^3x-y^4+2y^2x-y^3-2x^3+x^2y+3\in {\Bbb R}[x,y].$$
In Example \ref{Ej-asint-cuartica}, we obtain the  asymptote $\tilde{\cal D}$  of $\cal C$ at the infinity branch $B$ associated to $(1:0:0)$. We get that $\tilde{\cal D}$ is  defined by the polynomial
$$\tilde{f}(x,y)=y^3-x^2+y^2+1/3y+1/27\in {\Bbb R}[x,y].$$
From Proposition \ref{T-regular-asintota}, we deduce that $\tilde{\cal D}$ is a regular perfect curve, and all
%Observe that $P=(1:0:0)$ is a regular infinity point of $\tilde{\cal D}$ ($f_{d-1}(1, 0)\not=0$) (see proof of
%Corollary \ref{C-same-character}) and then  $\tilde{\cal D}$ is a regular perfect curve.
the asymptotes of $\cal C$ at $B$
are the curves defined by polynomials
$$y^3-x^2+y^2+ax+by+c,\,\,\,\,a,b,c\in {\Bbb R}.$$
These curves determine a proximity class isomorphic to $\mathbb{R}^{3}$ (Proposition \ref{P-evectorial}).

\para

In Figure \ref{F-inf-asymptotes}, we plot  $\cal C$, and the asymptotes $\tilde{\cal D}_i,\,i=1,2,3$ defined by the polynomials
$$\tilde{f}_1(x,y)=y^3-x^2+y^2-y+1/27,\qquad \tilde{f}_2(x,y)=y^3-x^2+y^2+1/3y-2x+1/27$$
$$\tilde{f}_3(x,y)=y^3-x^2+y^2.$$
%In Figure \ref{F-inf-asymptotes}, we plot  $\cal C$, and $\tilde{\cal D}_i,\,i=1,2,3$.
 \vspace*{-1.2cm}
\begin{figure}[h]\label{F-inf-asymptotes}
$$
\begin{array}{lcr}
\psfig{figure=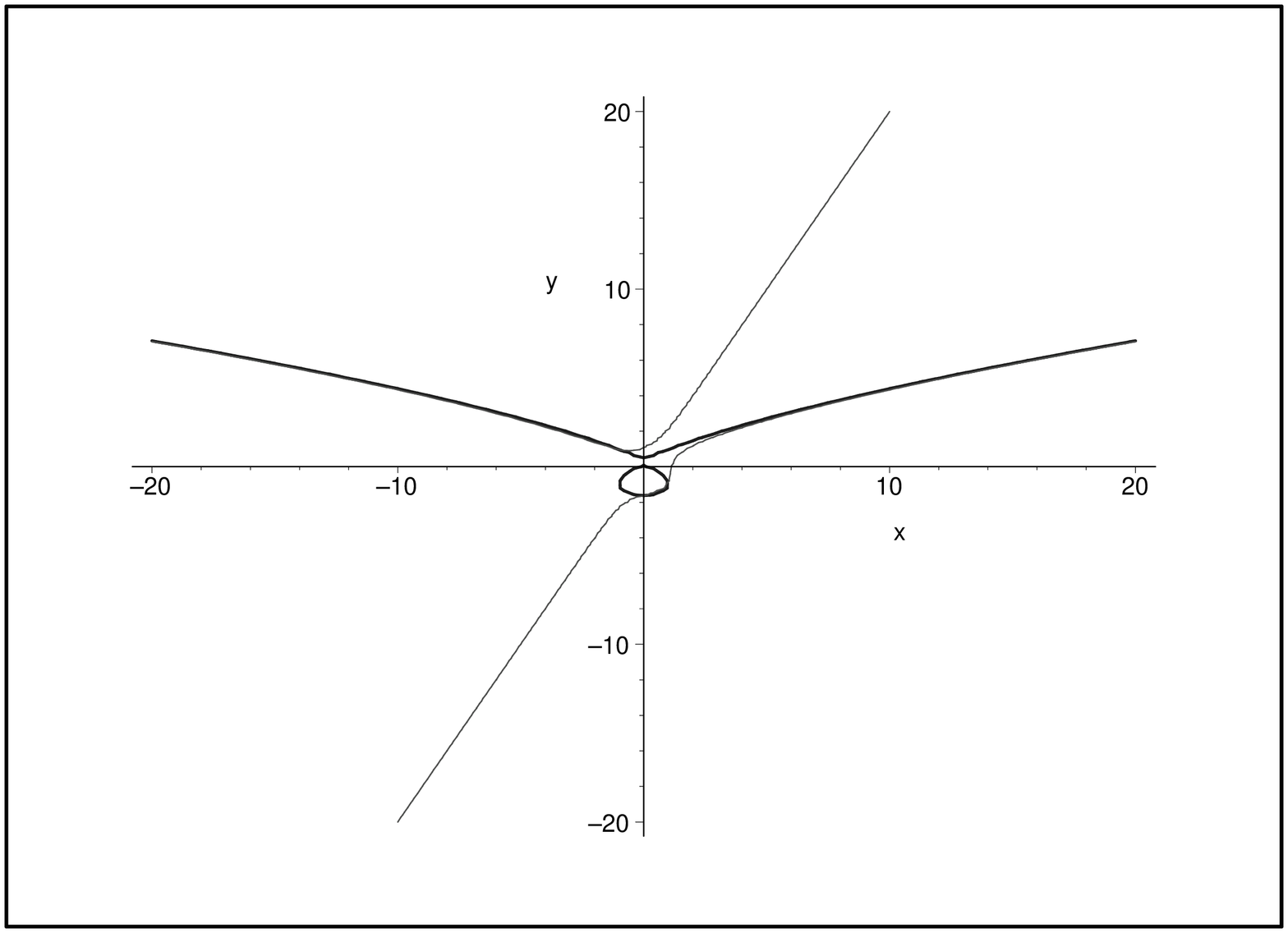,width=4.3cm,height=4.3cm,angle=270} &
\psfig{figure=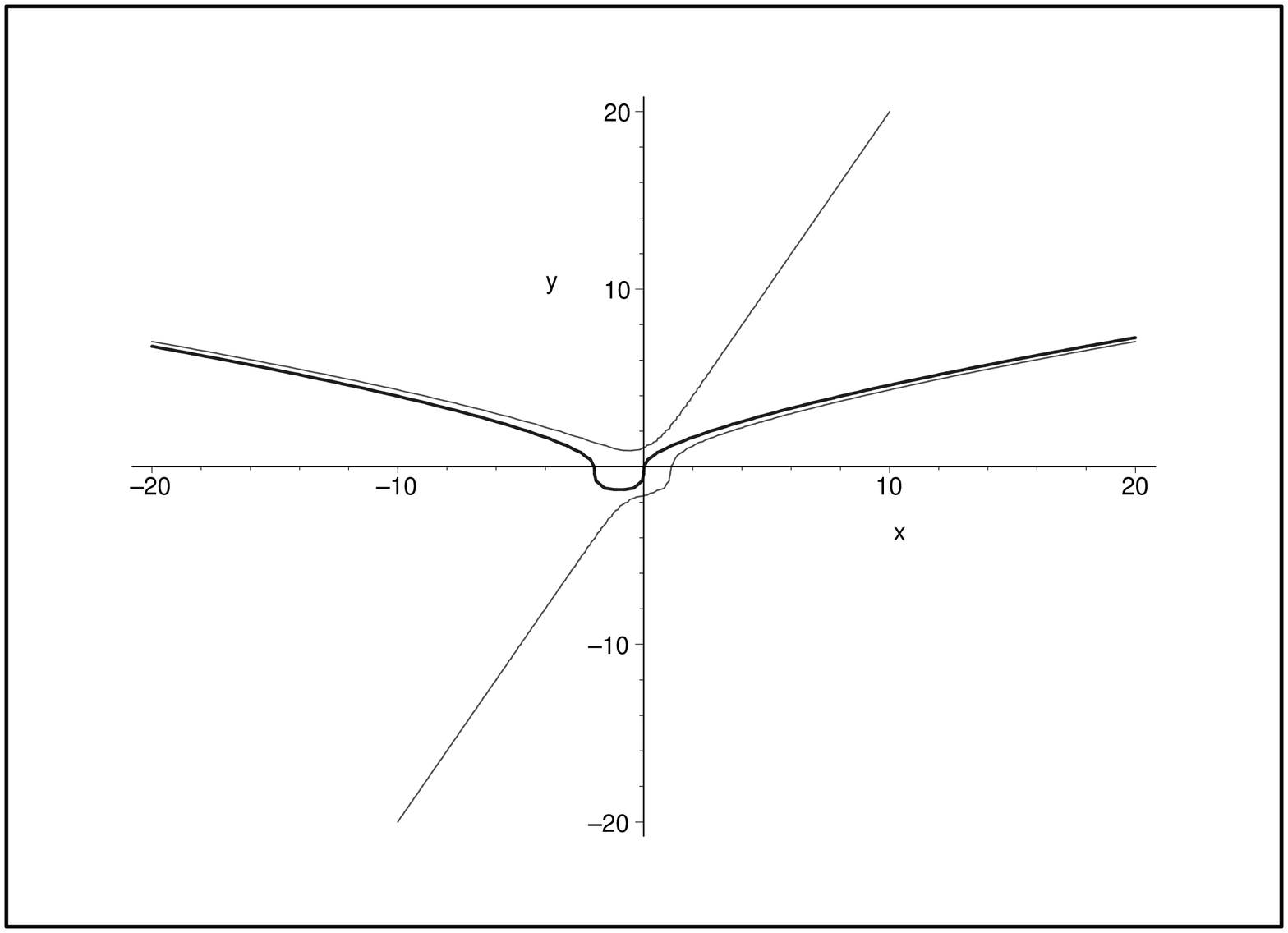,width=4.3cm,height=4.3cm,angle=270}&
\psfig{figure=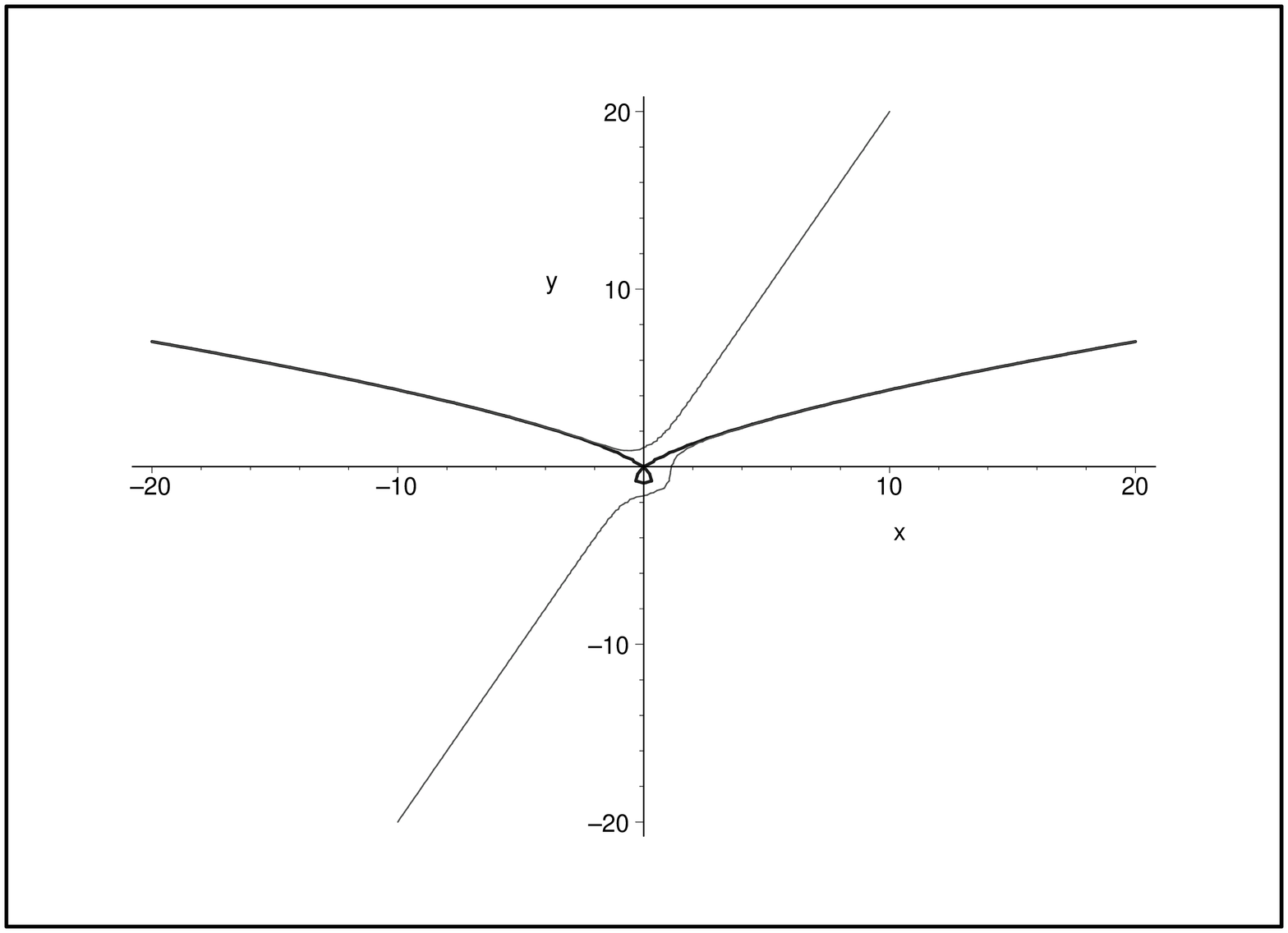,width=4.3cm,height=4.3cm,angle=270}
\end{array}
 \vspace*{-0.4cm}
$$ \caption{Curve $\cal C$, and asymptotes $\tilde{\cal D}_1$ (left), $\tilde{\cal D}_2$ (center) and $\tilde{\cal D}_3$ (right).}
\end{figure}
\end{example}

\end{document}